\documentclass[notitlepage,11pt,draft]{article}
\usepackage{amssymb,amsfonts,amsmath,geometry,esint,comment,upgreek,cases,graphicx,bbm,ulem}
\usepackage{scalerel,stackengine,mathtools,stmaryrd}
\usepackage{bm}

\usepackage{imakeidx}
\makeindex

\usepackage{color}

\catcode`\@=11
\@addtoreset{equation}{section}

\catcode`\@=12

\allowdisplaybreaks

\def\eps{\varepsilon}

\def\R{\mathbb{R}}
\def\f{\varphi}
\def\S{\mathbb S}
\def\rmH{{\rm H}}

\def\ird{\int\limits_{\R^d}}

\def\Ps{|\Pi\sigma|}

\def\div{{\rm div}}

\def\proof{\noindent{\textbf{Proof. }}}
\def\QED{\hfill {$\square$}\goodbreak \medskip}

\newtheorem{Theorem}{Theorem}[section]
\newtheorem{Lemma}[Theorem]{Lemma}
\newtheorem{Proposition}[Theorem]{Proposition}
\newtheorem{Corollary}[Theorem]{Corollary}
\newtheorem{Remark}[Theorem]{Remark}

\begin{document}

\title{
Hardy type inequalities with mixed cylindrical-spherical weights: \\
the general case}

\author{
Roberta Musina\footnote{Dipartimento di Scienze Matematiche, Informatiche e Fisiche, Universit\`a di Udine, Italy.
Email: {roberta.musina@uniud.it}, orcid.org/0000-0003-4835-8004. Partially supported by European Union -  
NextGenerationEU - PNRR M4.C2.1.1 - PRIN 2022 - 20227HX33Z\_003  - CUP 
G53D23001720006},
\setcounter{footnote}{6}
Alexander I. Nazarov\footnote{
St. Petersburg Department of Steklov Institute and St{.} Petersburg State University, St{.} Petersburg, Russia. E-mail: al.il.nazarov@gmail.com, 
orcid.org/0000-0001-9174-7000.}
}

\date{}

\maketitle

\bigskip

\noindent
{\small {\bf Abstract:} We continue our investigation of Hardy-type inequalities involving combinations of  cylindrical and spherical weights. 
Compared to [Cora-Musina-Nazarov, Ann. Sc. Norm. Sup., 2024], where the {\it quasi-spherical case} was considered, we handle the full range of allowed parameters. This has led to the observation of new phenomena related to lack of compactness.}



\section{Introduction and main results} 
We take integers $d>k\ge 1$ and put
$$
\R^d\equiv\R^k\times \R^{d-k}\ni (y,x)\equiv z~\!.
$$
Given $p>1$ and $a,b, \gamma\in\R$, we deal with the best constant $S_{b,\gamma}=S_{b,\gamma}(p,a)$ in dilation-invariant inequalities of the type
\begin{equation}
\label{eq:ine_MIX}
S_{b,\gamma} \int\limits_{\R^d}|y|^{a-p-b+\gamma}|z|^{-\gamma}|u|^p~\!dz\le \int\limits_{\R^d}|y|^{a}|z|^{-b}|\nabla u|^p~\!dz~,
\quad u\in {\cal C}^\infty_c(\R^d),
\end{equation}
under local summability hypotheses on the weights involved. That is, we assume that the following conditions
hold:
\begin{equation}
\label{eq:assumptions}
k+a>0~,\qquad d+a-p>b~,\qquad k+a-p>b-\gamma.
\end{equation}
Our starting motivation was born from the growing interest
in differential operators of the form
$$
-\div(|y|^aA(z)\nabla u)~\!,
$$
which are related to the the right hand side in (\ref{eq:ine_MIX}) for $p=2$ and for $A(z)$ controlled by $|z|^{-b}~\!\mathcal{I}_d$. 
Among the recent literature we cite \cite{CS, DFM, DM, STV1, STV2} and references therein.

The present paper is in fact a continuation of \cite{CMN}, where (\ref{eq:ine_MIX}) and similar inequalitiles in cones are investigated in the {\it quasi-spherical case}
$\gamma=p+b$. In particular, it has been proved there that 
\begin{equation}
\label{eq:CMN}
S_{b,p+b}=\inf_{u\in \mathcal D^{1,p}(\R^d;|y|^{a}|z|^{-b}~\!dz)\atop u=u(|z|)}
\frac{~\displaystyle\int\limits_{\R^d}|y|^{a}|z|^{-b}|\partial_r u|^p~\!dz}
{\displaystyle\int\limits_{\R^d}|y|^{a}|z|^{-b-p}|u|^p~\!dz}={\rmH_b^p}~,\qquad 
\text{where}~~\rmH_b=\frac{d+a-p-b}{p}>0~\!,
\end{equation}
provided that $k+a>0$ and $d+a-p>b$ (notice that the case of purely spherical weights is included by choosing $a=0$).
Here $\partial_ru=\nabla u\cdot z|z|^{-1}$ denotes the radial derivative of 
$u$, and the homogeneous weighted Sobolev-type space 
$$\mathcal D^{1,p}(\R^d;|y|^{a}|z|^{-b}~\!dz)$$ is the completion of ${\cal C}^\infty_c(\R^d)$ 
with respect to the norm
$$
\|u\|=\Big(\int\limits_{\R^d}|y|^{a}|z|^{-b}|\nabla u|^p dz\Big)^\frac1p.
$$
By \cite[Theorem 4]{CMN}, the infimum $S_{b,p+b}$ is not achieved on $\mathcal D^{1,p}(\R^d;|y|^{a}|z|^{-b}~\!dz)$.

\medskip

The value of $S_{b,\gamma}$ is also known in the {\it purely cylindrical case}: if $k+a>p$ we can take 
$b=\gamma=0$, and the Hardy inequality in $\R^k$ easily gives
\begin{equation}
\label{eq:cyl}
S_{0,0}= \inf_{u\in \mathcal D^{1,p}(\R^k;|y|^{a}dy)}
\frac{\displaystyle~\int\limits_{\R^k}|y|^a|\nabla u|^p~\!dy}
{\displaystyle\int\limits_{\R^k}|y|^{a-p}|u|^p ~\!dy}= \Lambda_0^p~,\qquad 
\text{where}~~\Lambda_0=\frac{k+a-p}{p} >0. 
\end{equation}
It is well known that $S_{0,0}$ is not achieved on $\mathcal D^{1,p}(\R^d;|y|^{a}~\!dz)$.

Our first result provides a necessary and sufficient condition to have that the best constant
$$
S_{b,\gamma}= \inf_{u\in \mathcal D^{1,p}(\R^d;|y|^{a}|z|^{-b} dz)}{\cal J}_{b,\gamma}(u)
$$
is positive, where
$$
{\cal J}_{b,\gamma}(u)= \frac{\displaystyle\int\limits_{\R^d}|y|^{a}|z|^{-b}|\nabla u|^p~\!dz}
{\displaystyle\int\limits_{\R^d}|y|^{a-p-b+\gamma}|z|^{-\gamma}|u|^p~\!dz}~\!.
$$

\begin{Theorem} 
\label{T:main_H}
Let (\ref{eq:assumptions}) holds. Then $S_{b,\gamma}>0$ if and only if 
$\gamma\ge b$.
\end{Theorem}
The hypothesis $\gamma\ge b$
 prevents vanishing 
along the singular set $\Sigma_0:=\{y=0\}$, a  phenomenon not previously observed in similar 
Hardy-type inequalities. 

Further information about the infimum $S_{b,\gamma}$ can be obtained by 
investigating the next minimization problem on the sphere, which deserves independent interest:
\begin{equation}
\label{eq:min_sphere}
\widetilde{S_{b,\gamma}}=
\inf_{\f \in W^{1,p}(\S^{d-1};\Ps^ad\sigma)} \widetilde{{\cal J}_{b,\gamma}}(\f)~\!.
\end{equation}
Here $\Pi:\R^d\to \R^k\times\{0\}$ is the orthogonal projection, 
$W^{1,p}(\S^{d-1};\Ps^ad\sigma)$ is a suitably defined non-homogeneous, weighted 
Sobolev space,
see Subsection \ref{S:sphere}, and
\begin{equation}
\label{eq:J_new}
\widetilde{{\cal J}_{b,\gamma}}(\f)= \frac{\displaystyle~\int\limits_{\S^{d-1}}|\Pi\sigma|^{a}\big[|\nabla_{\!\sigma} \f |^2+{\rmH_b^2}\f^2\big]^\frac{p}{2}~\!d\sigma~}
{\displaystyle\int\limits_{\S^{d-1}}|\Pi\sigma|^{a-p-b+\gamma}|\f |^p~\!d\sigma}~\!.
\end{equation}

Our results on the infima ${S_{b,\gamma}}, \widetilde{S_{b,\gamma}}$ are summarized in 
the next theorems. It is convenient to keep the case $\gamma>b$ separate from the {\it bottom case} $\gamma=b$, which includes the purely cylindrical case.

\begin{Theorem} 
\label{T:Main2}
Assume (\ref{eq:assumptions}) and let $\gamma>b$. Then 
$S_{b,\gamma}=\widetilde{S_{b,\gamma}}$. Moreover, 
the function $\gamma\mapsto S_{b,\gamma}$ is continuous 
and strictly increasing.
\end{Theorem}

The situation is more complex in the bottom case. 
Note that if $\gamma=b$ then the assumptions in (\ref{eq:assumptions}) transform into $k+a>p$ and $b<p\rmH_0
=d+a-p$, compare with  (\ref{eq:CMN}). Moreover, we have 
\begin{equation}
\label{eq:Sbb}
S_{b,b}= \inf_{u\in \mathcal D^{1,p}(\R^d;|y|^{a}|z|^{-b} dz)}
\frac{\displaystyle~\int\limits_{\R^d}|y|^{a}|z|^{-b}|\nabla u|^p~\!dz}
{\displaystyle\int\limits_{\R^d}|y|^{a-p}|z|^{-b}|u|^p~\!dz}~\!.
\end{equation}

\begin{Theorem}
\label{T:Main3}
Assume $k+a>p$. Then $S_{b,b}=\widetilde{S_{b,b}}$ for any $b<p\rmH_0$, and the 
following facts hold:
\begin{itemize}
\item[$i)$] The function $b\mapsto S_{b,b}$ is (locally) Lipschitz
continuous and non-increasing on $(-\infty,p\rmH_0)$; 
\item[$ii)$] There exists $b_*\in[0,p\rmH_0)$ such that
$$
\begin{cases}
~S_{b,b}=\Lambda_0^p&\text{if $b\le b_*$,}\\
~S_{b,b}<\Lambda_0^p&\text{if $b_*<b<p\rmH_0$,}
\end{cases}
$$
where $\Lambda_0^p$ is the Hardy constant in case $b=0$, see (\ref{eq:cyl});
\item[$iii)$] The function
$b\mapsto S_{b,b}$ is strictly decreasing for $b\in( b_*,p\rmH_0)$.
\end{itemize}
\end{Theorem}

In general, it is quite difficult to compute the 
explicit value of the infimum $S_{b,\gamma}$. This produces some technical (but challenging) difficulties in proving
the next result.

\begin{Theorem} 
\label{T:non-achiev}
Assume (\ref{eq:assumptions}) and let $\gamma\ge b$. Then $S_{b,\gamma}$ is not achieved on $\mathcal D^{1,p}(\R^d;|y|^{a}|z|^{-b}~\!dz)$.
\end{Theorem}

The paper is organized as follows. 
Section \ref{S:computations} contains some preliminaries and the
compactness Lemma \ref{L:compact}, which is crucially used to study the best constant $\widetilde{S_{b,\gamma}}$.

Section \ref{S:EL} includes some results of independent interest about positive (super)solutions to 
\begin{equation}
\label{eq:eq_m}
-\div(|y|^a|z|^{-b}|\nabla U|^{p-2}\nabla U)= 
m|y|^{a-p-b+\gamma}|z|^{-\gamma}|U|^{p-2}U ~\!.
\end{equation}

\medskip

We collect the proofs of the main results in Section \ref{S:proofs}, starting with 
Theorem \ref{T:main_H} in Subsection \ref{SS:main_H}. 
The proofs of Theorems \ref{T:Main2} and \ref{T:Main3} are somehow interconnected. Our strategy goes as follows:
\begin{enumerate}
\item Preliminary results on the infimum $\widetilde{S_{b,\gamma}}$ are given in Lemma \ref{L:main_sphere} of Subsection \ref{S:sphere}.
Thanks to Lemma \ref{L:compact}, the minimization problem (\ref{eq:min_sphere}) is compact if $\gamma>b$, so that $\widetilde{S_{b,\gamma}}$ is achieved. In the bottom case $\gamma=b$
remarkable phenomena related to lack of compactness are observed,
see 
Lemma \ref{L:widetilde_achieved};
\item 
The proof of Theorem \ref{T:Main2} in Subsection \ref{SS:Main2} is based on the previously proved results on $\widetilde{S_{b,\gamma}}$.
On the other hand, it also provides further crucial information on the infimum $\widetilde{S_{b,\gamma}}$ (some notable facts about the minimization problem in (\ref{eq:min_sphere}) are summarized in Remark \ref{R:final_tilde});
\item Finally, we prove 
Theorem \ref{T:Main3} in Subsection \ref{SS:Main3} and Theorem \ref{T:non-achiev} in Subsection \ref{SS:na}.
\end{enumerate}

\medskip

Section \ref{S:b*} contains estimates from below and from above on the borderline exponent $b_*$ appearing in the bottom case $\gamma=b$,
compare with Theorem \ref{T:Main3}. In the Hilbertian case $p=2$ the sharp value of $b_*$ and  the explicit value of $S_{b,b}$ are computed in 
Theorem \ref{T:sharp}.

\bigskip
 {\small
\noindent
{\bf Notation.}
The {\it singular set} $\Sigma_0=\{y=0\}$ is the kernel  of the orthogonal 
projection   
$\Pi(y,x)=(y,0)$; if $r>0, \sigma\in\S^{d-1}$ are the spherical coordinates of 
$z=(y,x)\in \R^k\times\R^{d-k}$, then
$|y|=r|\Pi\sigma|$.

We denote by $\mathcal K$  one of the 
following cones:
\begin{equation}
\label{eq:cone}
\begin{gathered}
\mathcal K=\R^{d}\setminus\Sigma_0, \quad
\text{ if $k+a\ge p$;}\\
\mathcal K=\R^{d}\setminus\{0\}, \quad \text{ if $k+a<p$.}
\end{gathered}
\end{equation}
Under the assumptions in (\ref{eq:assumptions}), for $\Omega\subseteq \mathcal K$ we put 
$$
S_{b,\gamma}(\Omega):= \inf_{u\in {\cal C}^\infty_c(\Omega)}
\frac{\displaystyle\int\limits_{\Omega}|y|^{a}|z|^{-b}|\nabla u|^p~\!dz}
{\displaystyle\int\limits_{\Omega}|y|^{a-p-b+\gamma}|z|^{-\gamma}|u|^p~\!dz}~\!.
$$

\medskip
By 
$\nabla_{\!\sigma}, \Delta_\sigma$ we denote the gradient and the Laplace-Beltrami operator on $\S^{d-1}$. 

\medskip

Let ${\mathfrak m}$ be a nonnegative and locally integrable weight in the open set $\Omega\subseteq\R^d$. We adopt standard notation for the weighted 
space $L^p(\Omega;{\mathfrak m}(z)dz)$, which is endowed with the norm
$$
 \|u\|_{p}:=\Big(\int\limits_{\Omega}|u|^p~\!{\mathfrak m}(z)dz\Big)^\frac1p.
$$
For $u\in {\cal C}^\infty(\Omega)$ we put
$$
\|u\|_{1,p}:= \|\nabla u \|_{p}+\|u\|_{p},
$$
and let  $W^{1,p}(\Omega;{\mathfrak m}(z)dz)$ be the completion of the set
$\big\{u\in {\cal C}^\infty(\Omega)~|~\|u\|_{1,p}<\infty \big\}$
with respect to the norm
$\|\cdot \|_{1,p}$. 
Next, $W^{1,p}_0(\Omega;{\mathfrak m}(z)dz)$ is the closure of ${\cal C}^\infty_c(\Omega)$ in 
$W^{1,p}(\Omega;{\mathfrak m}(z)dz)$. The space $W^{1,p}_{\rm loc}(\Omega;{\mathfrak m}(z)~\!dz)$ is introduced in a standard way.
If ${\mathfrak m}\equiv 1$ we suppress the volume element $dz$. 

Similar notation is used for weighted $L^p$ and non-homogeneous $W^{1,p}$-type spaces on the sphere. 

\medskip

We warn the reader that this notation differs from the standard one. For example, in \cite[Section 1.9]{HKM} the notation $H^{1,p}(\Omega;{\mathfrak m}(z)dz)$, 
$H^{1,p}_0(\Omega;{\mathfrak m}(z)dz)$ is used.

\medskip

We point out, see for 
instance \cite[Lemma 1]{CMN}, that for $\alpha, \beta\in \R$ we have 
$$
\begin{aligned}
&\text{$|y|^\alpha|z|^{-\beta}\in L^1_{\rm loc}(\R^d)$}&& \text{if and only if 
$k+\alpha>0$ and $d+\alpha>\beta$;}\\
&\text{$|\Pi \sigma|^\alpha\in L^1(\S^{d-1})$}&& \text{if and only if $k+\alpha>0$.}
\end{aligned}
$$

Through the paper, any positive constant whose value is not important is denoted 
by $c$. It may take different values at different places. To indicate that a 
constant depends on some parameters 
we list them in parentheses. 

}

\section{Preliminaries}
\label{S:computations}

\begin{Lemma}
\label{L:density}
Assume that $k+a>0$, $d+a>p+b$ and let $u\in {\cal C}^\infty_c(\R^d)$. Then 
there exists a sequence 
$u_h\in {\cal C}^\infty_c(\R^d\setminus\{0\})$ such that
$$
u_h\to u\quad \text{in} \quad \mathcal D^{1,p}(\R^d;|y|^a|z|^{-b}dz); \qquad
u_h\to u\quad \text{in} \quad L^p(\R^d;\mathfrak m(z)dz)
$$
for any  weight 
$\mathfrak m\in L^1_{\rm loc}(\R^d)$. If in addition $k+a\ge p$, then the 
sequence $u_h$ 
can be taken in ${\cal C}^\infty_c(\R^d\setminus\Sigma_0)$.
\end{Lemma}

\proof
Take $\eta\in  {\cal C}^\infty(\R)$ such that $\eta\equiv 1$ on $(-\infty,1]$, 
$\eta\equiv 0$ on $[2,\infty)$ and $|\eta|\le 1$. For any integer $h\ge 
1$ 
put
$$
u_h(y,x)=
\begin{cases}
\eta\Big(\dfrac{-\log |z|}{h}\Big)u(y,x)&\text{if \ $0<k+a<p$;}\\ \\
\eta\Big(\dfrac{-\log |y|}{h}\Big)u(y,x)&\text{if \ $k+a\ge p$.}
\end{cases}
$$
Then $u_h\in  {\cal C}_c^\infty(\R^d\setminus\{0\})$ and $u_h\in  {\cal C}_c^\infty(\R^d\setminus\Sigma_0)$ if $k+a\ge p$.

Evidently, $u_h\to u$ pointwise. Moreover, by standard computations (cf. the 
proof of Theorem 2 in \cite{CMN}) one can prove that 
$$
 \begin{aligned}
 \int\limits_{\R^d}|u_h|^p\mathfrak 
m(z)~\!dz&=\int\limits_{\R^d}|u|^p\mathfrak m(z)~\!dz+o(1),\\
 \int\limits_{\R^d}|y|^{a}|z|^{-{{b}}}|\nabla 
u_h|^p~\!dz&=\int\limits_{\R^d}|y|^{a}|z|^{-{{b}}}|\nabla u|^p~\!dz+o(1)
 \end{aligned}
$$
as $h\to\infty$.  This concludes the proof.
\QED

As a consequence of Lemma \ref{L:density} we obtain the following result.

\begin{Corollary}
\label{C:equal}
Under the assumptions in (\ref{eq:assumptions}) it holds that
$$
S_{b,\gamma}(\R^{d}\setminus\{0\})= S_{b,\gamma}.
$$
If in addition $k+a\ge p$, then
$$
S_{b,\gamma}(\R^{d}\setminus\Sigma_0)=S_{b,\gamma}.
$$
\end{Corollary}

Next, we deal with functions on the sphere.
Let $\{{\rm e}_j\}$ be the canonical basis of $\R^d$ and let 
$\sigma_j=\sigma\cdot {\rm e}_j$ be the $j$-th coordinate of $\sigma 
\in\S^{d-1}$ in $\R^d$. Then
$$
\nabla_{\!\sigma}\sigma_j={\rm e}_j-\sigma_j\sigma~\!,\quad 
|\nabla_{\!\sigma}\sigma_j|^2=1-\sigma_j^2~\!.
$$
Since $\Ps^2=\sum\limits_{j=1}^k\sigma_j^2$, it follows that $\nabla_{\!\sigma}\Ps=\Ps^{-1}\Pi\sigma-\Ps\sigma$ on $\S^{d-1}\setminus\Sigma_0$, hence
\begin{equation}
\label{eq:useful}
|\nabla_{\!\sigma}\Ps|^2=1-\Ps^2\qquad \text{on $\S^{d-1}\setminus\Sigma_0$.}
\end{equation}

\medskip

\begin{Lemma}
\label{L:useful_new}
Assume $k+a>0$. 
Then $\Ps^{-\lambda}\in W^{1,p}(\S^{d-1};\Ps^ad\sigma)$ for any 
$\lambda<\Lambda_0$, where $\Lambda_0$ is introduced in (\ref{eq:cyl}).
\end{Lemma}

\proof 
Recall that $\Ps^\tau\in L^1(\S^{d-1})$ if and only if $k+\tau>0$. 
Together with (\ref{eq:useful}), this gives
$$
\Ps^{-\lambda}, |\nabla_{\!\sigma}\Ps^{-\lambda}|\in L^p(\S^{d-1};\Ps^ad\sigma)\quad
\text{if and only if $\lambda<\Lambda_0$}.
$$
To conclude the proof notice that the smooth functions 
$\f_\eps(\sigma)=(\eps^2+\Ps^2)^{-\frac\lambda2}$ approximate
$\Ps^{-\lambda}$ in $W^{1,p}(\S^{d-1};\Ps^ad\sigma)$ as $\eps\to 0$.
\QED

\begin{Lemma}
\label{L:useful_iii}
If $k+\tau>2$, then
$$
(k+\tau-2) \int\limits_{\S^{d-1}}|\Pi\sigma|^{\tau-2}  d\sigma 
=(d+\tau-2)\int\limits_{\S^{d-1}}  |\Pi\sigma|^{\tau}d\sigma.
$$
\end{Lemma}

\proof 
As before,  let $\sigma_j$ be the $j$-th coordinate of $\sigma \in\S^{d-1}\subset \R^d$.
It is well known that $\sigma_j$ is an eigenfunction for the Laplace-Beltrami 
operator on $\S^{d-1}$ relative
to the eigenvalue $(d-1)$, so we can compute
$$
-\frac12\Delta_\sigma\sigma_j^2=-\sigma_j\Delta_\sigma\sigma_j-|\nabla_{\!\sigma
}
\sigma_j|^2=(d-1)\sigma_j^2-(1-\sigma_j^2).
$$
It follows that
$$
-\frac12\Delta_\sigma\Ps^2= -\frac12\sum\limits_{j=1}^k\Delta_\sigma\sigma_j^2= 
d\sum\limits_{j=1}^k\sigma_j^2-k=d\Ps^2-k.
$$
Therefore, integration by parts gives
$$
\begin{aligned}
k\int\limits_{\S^{d-1}} |\Pi\sigma|^{\tau-2}  d\sigma&=
d\int\limits_{\S^{d-1}} |\Pi\sigma|^{\tau}  d\sigma+ \frac12 
\int\limits_{\S^{d-1}}\Ps^{\tau-2}\Delta_\sigma\Ps^2~\!d\sigma\\
&= d\int\limits_{\S^{d-1}} |\Pi\sigma|^{\tau}  d\sigma-(\tau-2) 
\int\limits_{\S^{d-1}} \Ps^{\tau-2}|\nabla_{\!\sigma}\Ps|^2  d\sigma
\\&
= d\int\limits_{\S^{d-1}} |\Pi\sigma|^{\tau}  d\sigma-(\tau-2) 
\int\limits_{\S^{d-1}} \Ps^{\tau-2}(1-\Ps^2)  d\sigma
\end{aligned}
$$
by (\ref{eq:useful}). The conclusion follows.
\QED

\begin{Lemma} 
\label{L:continuous}
Let $\omega$ be open in $\S^{d-1}$ 
and assume that
$\S^{d-1}\cap \Sigma_0\not\subset \overline\omega$.
\begin{itemize}
\item[$i)$] Assume $k+a>p$. Then $W^{1,p}_0(\omega;|\Pi\sigma|^ad\sigma)$ is continuously embedded into $L^p(\omega;|\Pi\sigma|^{a-p} d\sigma)$;

\item[$ii)$] Assume $k+a>0$ and let $k+\tau>0, \tau> a-p$. Then $W^{1,p}_0(\omega;|\Pi\sigma|^ad\sigma)$ is compactly embedded into
$L^p(\omega;|\Pi\sigma|^\tau d\sigma)$.
\end{itemize}
 \end{Lemma}

\proof
We can assume that the north pole  ${\rm e}_d=(0,\dots 0,1)$ does 
not belong to 
$\overline\omega$.

Let ${\rm P}:\R^{d-1}\equiv  \R^k\times \R^{d-k-1}\to \S^{d-1}\setminus\{{\rm 
e}_d\}\subset \R^{d}$ be the inverse of the stereographic projection from 
${\rm e}_d$.
More explicitly, 
$$
{\rm P}(y,\xi)=(\mu y; \mu \xi, 1-\mu)\in \S^{d-1}\subset 
\R^k\times(\R^{d-k-1}\times \R)~,\quad 
\mu\equiv\mu(y,\xi)=\frac{2}{1+|y|^2+|\xi|^2}
$$
(the variable $\xi$ has to be omitted if $d=k+1$). 

For $\f\in {\cal C}^\infty_c(\omega)$ we put 
$\tilde{\f}:=\f\circ {\rm P}\in  {\cal C}^\infty_c(\Omega)$, where
$\Omega\subset\R^{d-1}$ is the (open and 
bounded) set such that ${\rm P}(\Omega)=\omega$. 
Since 
$\mu$ is bounded and bounded away from $0$ on $\Omega$, 
we have that there exists $c>1$ depending only on $\omega$ such that 
\begin{equation}
\label{eq:mu}
\begin{aligned}
\int\limits_{\omega}|\Pi \sigma|^a|\nabla_{\!\sigma} 
\f|^p~\!d\sigma&\,=
\int\limits_{\Omega}\mu^{d-1+a-p}|y|^a|\nabla {\tilde{\f}}|^p~\!d\xi dy\ge
c^{-1}\int\limits_{\Omega}|y|^a|\nabla {\tilde{\f}}|^p~\!d\xi dy
\\
\int\limits_{\omega\cap \{\Ps<\eps\}}|\Pi 
\sigma|^{\tau}|\f|^pd\sigma&\,\le
\int\limits_{\Omega\cap \{|y|<R\eps\} 
}\mu^{d-1+\tau}|y|^{\tau}|{\tilde{\f}}|^pd\xi dy
\le c \int\limits_{\Omega\cap \{|y|<R\eps\} }|y|^{\tau}|{\tilde{\f}}|^pd\xi dy
\end{aligned}
\end{equation}
for any $\eps>0$, provided that $k+\tau>0$, where $R>0$ is such that 
$\text{supp}\,\tilde\f\subset \Omega\subset B_R^{d-1}$. 

\medskip

To prove $i)$ we recall the Hardy inequality with cylindrical weights
$$
\Big(\frac{k+a-p}{p}\Big)^p \int\limits_{\R^{d-1}}|y|^{a-p}|{\tilde{\f}}|^pd\xi dy
\le \int\limits_{\R^{d-1}}|y|^{a}|\nabla {\tilde{\f}}|^p~\!d\xi dy,
$$
compare with (\ref{eq:cyl}). Taking $\tau=a-p$ and $\eps=1$ in (\ref{eq:mu}), we see that there exists $c>0$ depending on $R$ (hence on $\omega$), such that
$$
\displaystyle
\int\limits_{\omega}|\Pi\sigma|^{a-p}|\f|^p~\!d\sigma\le 
c\int\limits_{\omega}|\Pi\sigma|^a|\nabla_{\!\sigma}
\f|^p~\!d\sigma
$$
for any $\f \in {\cal C}^\infty_c(\omega)$. Statement $i)$ follows by recalling that  ${\cal C}^\infty_c(\omega)$ is dense in $W^{1,p}_0(\omega;|\Pi\sigma|^ad\sigma)$.

\medskip

To prove $ii)$ take a small $\delta>0$, in such a way that $k+\tau-\delta>0$ and $\tau-\delta+p>a$. We will now use the
Hardy inequality in the  form
$$
\Big(\frac{k+\tau-\delta}{p}\Big)^p \int\limits_{\R^{d-1}}|y|^{\tau-\delta}|{\tilde{\f}}|^pd\xi dy
\le \int\limits_{\R^{d-1}}|y|^{\tau-\delta+p}|\nabla {\tilde{\f}}|^p~\!d\xi dy~\!,
$$
to get
$$
\int\limits_{\Omega\cap \{|y|<R\eps\} }\!\!\!|y|^{\tau}|{\tilde{\f}}|^pd\xi dy\le 
c\eps^\delta \!\!\!\int\limits_{\Omega\cap \{|y|<R\eps\} }\!\!\!|y|^{\tau-\delta}|{\tilde{\f}}|^pd\xi dy\le 
c\eps^{\delta}\int\limits_{\Omega}|y|^{\tau-\delta+p}|{\tilde{\f}}|^pd\xi dy
\le c\eps^{\delta}\int\limits_{\Omega}|y|^{a}|\nabla {\tilde{\f}}|^pd\xi dy~\!.
$$
Together with (\ref{eq:mu}), this implies that
there exists 
 $c>0$  depending on $ \omega$, $\tau$ and $\delta$, such that 
\begin{equation}
\label{eq:ball}
\displaystyle
\int\limits_{\omega\cap \{\Ps<\eps\}}|\Pi\sigma|^{\tau}|\f|^p~\!d\sigma\le 
c\eps^{\delta}\int\limits_{\omega}|\Pi\sigma|^a|\nabla_{\!\sigma}
\f|^p~\!d\sigma
\end{equation}
for any $\eps>0$.
By density, we infer that (\ref{eq:ball}) holds for any $\f \in  W^{1,p}_0(\omega;|\Pi\sigma|^ad\sigma)$.

\medskip

To conclude the proof of $ii)$, we notice that for any $\eps>0$ the embedding operator 
$$W^{1,p}_0(\omega;|\Pi\sigma|^ad\sigma)\hookrightarrow 
L^p(\omega\cap \{\Ps>\eps\};|\Pi\sigma|^\tau d\sigma)
$$
is compact by the Rellich theorem. Therefore, the estimate 
(\ref{eq:ball}) shows that the embedding  operator
$W^{1,p}_0(\omega;|\Pi\sigma|^ad\sigma)\hookrightarrow 
L^p(\omega;|\Pi\sigma|^\tau d\sigma)$
can be approximated in norm by compact operators. Thus it is compact itself.
This completes the proof of the lemma.
\QED

\begin{Lemma}
\label{L:compact}
Let $k+a>0$ and $k+\tau>0$.
The space $W^{1,p}(\S^{d-1};|\Pi\sigma|^ad\sigma)$ is continuously embedded in 
$L^p(\S^{d-1};|\Pi\sigma|^\tau d\sigma)$
if and only if  $\tau\ge a-p$.  
This embedding is compact if $\tau>a-p$. 
\end{Lemma}

\proof
Assume $\tau\ge a-p$.
Take a point ${\rm e}_d\in\S^{d-1}\cap \Sigma_0$ and a cut-off function 
$\eta\in 
{\cal C}^\infty_c(\S^{d-1}\setminus\{{\rm e}_d\})$
such that $\eta\equiv 1$ in a neighborhood of $-{\rm e}_d$. By Lemma 
\ref{L:continuous}, the operators
$$
\f\mapsto \eta\f~,\quad \f\mapsto (1-\eta)\f~,\qquad  
W^{1,p}_0(\S^{d-1};|\Pi\sigma|^ad\sigma)\to 
L^p(\S^{d-1};|\Pi\sigma|^\tau d\sigma)
$$
are continuous if $\tau\ge a-p$, and compact if $\tau>a-p$. 

If $-k<\tau<a-p$ then the function $\Ps^{-\frac{k+\tau}{p}}$ belongs to 
$W^{1,p}(\S^{d-1};\Ps^ad\sigma)$ by Lemma \ref{L:useful_new} but evidently
$\Ps^{-\frac{k+\tau}{p}}\notin L^p(\S^{d-1};|\Pi\sigma|^\tau d\sigma)$. This 
completes the proof.
\QED

\section{On the Euler-Lagrange equations}
\label{S:EL}

Let $\mathcal K$ be the cone in (\ref{eq:cone}).
 We say that 
$U\in W^{1,p}_{\rm loc}({\mathcal K}; |y|^a|z|^{-b}~\!dz)$ is a  weak supersolution of the equation (\ref{eq:eq_m}), if 
\begin{equation}
\label{eq:EL_m}
\int\limits_{\mathcal K} |y|^a|z|^{-b}|\nabla U|^{p-2}\nabla U\cdot \nabla 
v~\!dz\ge m\int\limits_{\mathcal 
K}|y|^{a-p-b+\gamma}|z|^{-\gamma}|U|^{p-2}Uv~\!dz
\end{equation}
for any nonnegative function $v\in W^{1,p}({\mathcal K}; |y|^a|z|^{-b}~\!dz)$ 
compactly supported in $\mathcal K$. If equality holds in (\ref{eq:EL_m}) for any function $v\in W^{1,p}({\mathcal K}; |y|^a|z|^{-b}~\!dz)$, we say that 
$U$ is a weak solution  of the equation (\ref{eq:eq_m}).

\begin{Lemma}
\label{L:lsc}
Let $U\in W^{1,p}_{\rm loc}({\mathcal K}; |y|^a|z|^{-b}~\!dz)\setminus\{0\}$ be a nonnegative weak solution of (\ref{eq:eq_m}), for some 
$m\ge 0$. Then $U$ is positive and $\nabla U$ is (locally) H\"older continuous on $\R^d\setminus\Sigma_0$.
If in addition, if $k+a<p$,
then $U$ is lower semicontinuous and positive on $\R^d\setminus\{0\}$.
\end{Lemma}

\proof
On the cone $\R^d\setminus\Sigma_0$ the weights involved are locally bounded and bounded away from $0$. Thus 
$U$ is  of class ${\cal C}^{1,\alpha}$ locally on $\R^d\setminus\Sigma_0$
by the results in \cite{DiB, Tol} and  positive  thanks to the maximum principle. If in addition $k+a<p$,
then $|y|^a$ belongs to the Muckenhoupt class $A_p$. 
Since evidently 
\begin{equation*}
-{\rm div}\big(|y|^a|z|^{-b}\,|\nabla U|^{p-2}\nabla U\big)\ge 0\qquad \text{in 
$\R^d\setminus \{0\}$,}
\end{equation*}
we  conclude that $U$ is lower semicontinuous and positive in $\R^d\setminus\{0\}$ by \cite[Theorems 3.51, 3.66]{HKM}.
\QED

The following result includes adaptations of the Pinchover--Tintarev \cite{PT} and of the Allegretto--Piepenbrink 
\cite{Al,Pi} basic ideas. 
 
\begin{Theorem}
\label{T:AP}
Let the assumptions 
(\ref{eq:assumptions}) and $\gamma\ge b$ be fulfilled, and let 
$m> 0$. Then the following statements are equivalent:
\begin{itemize}

\item[$i)$] $S_{b,\gamma}({\mathcal K})\ge m$ (by Corollary  \ref{C:equal} this is 
equivalent to $S_{b,\gamma}\ge m$);

\item[$ii)$] There exists a lower semicontinuous
and positive weak solution $U\in W^{1,p}_{\rm loc}({\mathcal K}; |y|^a|z|^{-b}~\!dz)$ of equation (\ref{eq:eq_m});

 \item[$iii)$] There exists a lower semicontinuous
and positive weak supersolution $U\in W^{1,p}_{\rm loc}({\mathcal K}; |y|^a|z|^{-b}~\!dz)$ of equation (\ref{eq:eq_m}).

\end{itemize}
\end{Theorem}

\proof
The implication $ii)~\Rightarrow~iii)$ does not need a proof. 

To prove that $i)~\Rightarrow~ii)$ we first consider the case $k+a\ge p$ and $\mathcal 
K=\R^{d}\setminus\Sigma_0$. We introduce the sequence of expanding domains 
$$
\Omega_h=\{z\in\R^d\ :\ h^{-1}<|z|<2h,\ |y|>h^{-1}\}.
$$
It is evident that $S_{b,\gamma}(\Omega_h)> m$, and thus the form
\begin{equation}
\label{eq:norm}
{\cal E}_h^p[u]
:=\int\limits_{\Omega_h}|y|^a|z|^{-b}|\nabla u|^p~\!dz 
-m\int\limits_{\Omega_h}|y|^{a-p-b+\gamma}|z|^{-\gamma}|u|^p\,dz
\end{equation}
is equivalent to the $p$-th power of the norm in 
$W^{1,p}_0(\Omega_h;|y|^a|z|^{-b}dz)$. Therefore, 
for any nonnegative  $f_h\in{\cal 
C}^\infty_c(\Omega_h\setminus\Omega_{h-1})\setminus\{0\}$,
the functional 
$$
\frac1p~\!{\cal E}_h^p[u]-\int\limits_{\Omega_h}f_hu~\!dz
$$
has a unique (nonnegative) minimizer $u_h\in 
W^{1,p}_0(\Omega_h;|y|^a|z|^{-b}dz)$ that is a weak solution of
$$
-\div(|y|^a|z|^{-b}|\nabla u_h|^{p-2}\nabla u_h)- 
m|y|^{a-p-b+\gamma}|z|^{-\gamma}|u_h|^{p-2}u_h=f_h \quad \text{in 
$\Omega_h$},\qquad u_h|_{\partial \Omega_h}=0~\!.
$$
 By  elliptic regularity \cite{DiB, Tol}, we have that $u_h\in {\cal 
C}^{1,\alpha}(\Omega_h)$. 

Fix a point $z_0\in\Omega_1$. The Harnack inequality \cite{Ser, Tr} implies 
that $u_h>0$ in $\Omega_h$, and that the sequence 
$U_h(z)=u_h(z)/u_h(z_0)$ admits a subsequence which 
converges locally uniformly in $\mathcal K$ to a positive
weak solution $U$ of the equation 
(\ref{eq:eq_m}). Thus $i)~\Rightarrow~ii)$ is proved in this case.
\medskip

The case $k+a< p$ and $\mathcal 
K=\R^{d}\setminus\{0\}$ is more involved. Now we define 
$$
\Omega_h=\{z\in\R^d\ :\ h^{-1}<|z|<2h\}.
$$
Since the form (\ref{eq:norm}) is equivalent to the 
$p$-th power of the norm in $W^{1,p}_0(\Omega_h;|y|^a|z|^{-b}dz)$, we can
define $u_h\in W^{1,p}_0(\Omega_h;|y|^a|z|^{-b}dz)$ as  in previous case. Fix $z_0\in\Omega_1$.
The sequence $U_h(z)=u_h(z)/u_h(z_0)$
is locally bounded in $W^{1,p}(\R^d\setminus\{0\};|y|^a|z|^{-b}dz)$. Thus it
admits a subsequence which 
converges weakly in $W^{1,p}_{\rm loc}(\R^d\setminus\{0\};|y|^a|z|^{-b}dz)$ to 
a 
weak solution $U$ of the equation (\ref{eq:eq_m}). 
The Harnack inequality can be used to prove that $U_h$ is uniformly bounded away from $0$
on a compact subset of $\Omega_1$. Thus $U\neq 0$. Then Lemma \ref{L:lsc} gives
that $U$
is lower semicontinuous and positive in $\R^d\setminus\{0\}$. 
The implication $i)~\Rightarrow~ii)$ is proved also in this case.

\medskip

Finally, assume that $iii)$ holds. 
For any $u\in {\cal C}^\infty_c({\mathcal K})$ we use $v=U^{1-p}|u|^p$ as a
test function in (\ref{eq:EL_m}). Then we apply the 
Cauchy and Young 
inequalities to get
$$
\begin{aligned}
m\int\limits_{\R^d}&|y|^{a-p-b+\gamma}|z|^{-\gamma}|u|^p\,dz
\le \int\limits_{{\R^d}}|y|^a|z|^{-b}|\nabla U|^{p-2}\nabla 
U\cdot\nabla(U^{1-p}|u|^p) \,dz
\\&
\le \int\limits_{{\R^d}}|y|^a|z|^{-b}\Big(p\,\big|u U^{-1}\nabla 
U|^{p-1}|\nabla u|-(p-1)|u U^{-1}\nabla U|^p\Big)\,dz
\le \int\limits_{{\R^d}}|y|^a|z|^{-b}|\nabla u|^p
~\!dz,
\end{aligned}
$$
which proves the implication $iii)~\Rightarrow~i)$, thanks to the density result in Lemma \ref{L:density}.
\QED

\begin{Lemma}
\label{L:AP2}
Let the assumptions 
(\ref{eq:assumptions}) and $\gamma\ge b$ be fulfilled, and let $\mathcal K$ be the cone in (\ref{eq:cone}). 
Let $U\in W^{1,p}_{\rm loc}({\mathcal K}; |y|^a|z|^{-b}~\!dz)$ be a nonnegative
solution 
of the equation 
\begin{equation}
\label{eq:eq_mS}
-\div(|y|^a|z|^{-b}|\nabla v|^{p-2}\nabla v)=  S_{b,\gamma}
|y|^{a-p-b+\gamma}|z|^{-\gamma}v^{p-1}
\end{equation}  
in $\mathcal K$.
If $U\notin \mathcal D^{1,p}(\R^d;|y|^a|z|^{-b}dz)$, then $S_{b,\gamma}$ is not achieved in $\mathcal D^{1,p}(\R^d;|y|^a|z|^{-b}dz)$.
\end{Lemma}

\proof
We adapt the argument in \cite[Theorem 4]{CMN}. 
Assume by contradiction that the infimum $S_{b,\gamma}$ is achieved by
$u\in {\mathcal D}^{1,p}(\R^d;|y|^a|z|^{-b}dz)$. Up to a change of sign, the minimizer $u$ is  nonnegative
and solves (\ref{eq:eq_mS}) in a weak sense on $\R^d$. Then Lemma \ref{L:lsc} implies that the functions
$U$ and $u$ are positive on $\mathcal K$ and  of class ${\cal C}^1$ on 
$\R^d\setminus\Sigma_0$.

Next, take a domain $A$, compactly contained in $\R^d\setminus\Sigma_0$. As in 
the proof of \cite[Theorem 4]{CMN}, we  show that the open sets
$$W_1=\big\{z\in A~|~~ \nabla u\cdot\nabla U < |\nabla 
u||\nabla U| \big\}~,\quad
W_2=\Big\{z\in A~|~~
\frac{|\nabla u|}{u}\neq \frac{|\nabla U|}{U}~\Big\}
$$
are empty. By the arbitrariness of $A$, this easily imply that there exists a constant $\lambda>0$ such that 
$U=\lambda u$  on 
$\R^d\setminus\Sigma_0$, which contradicts the assumption $U\notin {\mathcal D}^{1,p}(\R^d;|y|^a|z|^{-b}dz)$ and concludes the proof.

\medskip

To show that the sets $W_1, W_2$ are empty we adapt the argument in \cite{CMN}.
Use Lemma \ref{L:density} to find a sequence of nonnegative functions $u_h\in C^\infty_c(\mathcal K)$ such that $u_h\to u$ in $\mathcal D^{1,p}(\R^d; |y|^a|z|^{-b}dz)$.
Since $U$ solves (\ref{eq:eq_mS}) in $\mathcal K$ and is strictly positive on the support of $u_h$, by using $U^{1-p}u_h^p$ as test function we have that
$$
S_{b,\gamma}
\ird|y|^{a-p-b+\gamma}|z|^{-\gamma}~\!u_h^p~\!dz\le  \ird|y|^{a}|z|^{-b}~|\nabla U|^{p-2}\nabla U\cdot \nabla \big(U^{1-p}\!u_h^p\big)~\!dz.
$$
Next, we use $u_h\to u$ in $L^p(\R^d;|y|^{a-p-b+\gamma}|z|^{-\gamma}dz)$ in the left-hand side and perform direct calculation to infer 
\begin{multline}
\label{eq:A1}
S_{b,\gamma}
\int\limits_{{\mathcal K}}|y|^{a-p-b+\gamma}|z|^{-\gamma}u^p~\!dz\\
=
\int\limits_{{\mathcal K}}|y|^a|z|^{-b}|\nabla U|^{p-2}
\Big(p~\!\frac{\nabla U\cdot\nabla u_h}{U^{p-1}} u_h^{p-1}
-(p-1)\frac{|\nabla U|^2u_h^p}{U^p}
\Big)dz+o(1)
\end{multline}
as $h\to \infty$. Next, we put
\begin{multline*}
\nu_1:=p\int\limits_{W_1}|y|^a|z|^{-b}|\nabla U|^{p-2}
\Big(\frac{|\nabla U|~\!|\nabla u|}{U^{p-1}}-\frac{\nabla U\cdot\nabla u}{U^{p-1}}\Big) ~\!u^{p-1}~\!dz\\=
p\int\limits_{W_1}|y|^a|z|^{-b}|\nabla U|^{p-2}
\frac{|\nabla U|~\!|\nabla u_h|-\nabla U\cdot\nabla u_h}{U^{p-1}} ~\!u_h^{p-1}~\!dz+o(1).
\end{multline*}
By the Cauchy inequality in $\R^d$, we have 
that the 
nonnegative constant $\nu_1$ is positive if and only if $W_1$ has positive 
measure. 
Thus, from (\ref{eq:A1}) it follows that 
\begin{multline}
\label{eq:Young}
S_{b,\gamma}
\int\limits_{{\mathcal K}}|y|^{a-p-b+\gamma}|z|^{-\gamma}u^p~\!dz\\
\le 
\int\limits_{{\mathcal K}}|y|^a|z|^{-b}
\bigg[p\Big(\frac{|\nabla U|u_h}{U}\Big)^{p-1}|\nabla u_h|
-(p-1)
\Big(\frac{|\nabla U|u_h}{U}\Big)^{p}\bigg]dz - \nu_1+o(1).
\end{multline}
Further, as $h\to \infty$ we have 
\begin{multline*}
\nu_2:=
\int\limits_{W_2}|y|^a|z|^{-b}
\Big(|\nabla u|^p+(p-1)\Big(\frac{|\nabla U|u}{U}\Big)^p-
p\Big(\frac{|\nabla U|u}{U}\Big)^{p-1}|\nabla u|\Big)dz\\
=\int\limits_{W_2}|y|^a|z|^{-b}
\Big(|\nabla u_h|^p+(p-1)\Big(\frac{|\nabla U|u_h}{U}\Big)^p-
p\Big(\frac{|\nabla U|u_h}{U}\Big)^{p-1}|\nabla u_h|\Big)dz+o(1)~\!.
\end{multline*}
By Young's inequality, it turns out that $\nu_2>0$ if and only if $W_2$ has positive measure.
Thus, from (\ref{eq:Young}) it follows that 
$$
S_{b,\gamma}
\int\limits_{{\mathcal K}}|y|^{a-p-b+\gamma}|z|^{-\gamma}~\!u^p~\!dz
\le 
\int\limits_{{\mathcal K}}|y|^a|z|^{-b}|\nabla u_h|^p~\!dz - (\nu_1+\nu_2)+o(1).
$$
Letting $h\to \infty$ we obtain
$$
\begin{aligned}
S_{b,\gamma}
\int\limits_{{\mathcal K}}|y|^{a-p-b+\gamma}|z|^{-\gamma}~\!u^p~\!dz
&\le 
\int\limits_{{\mathcal K}}|y|^a|z|^{-b}|\nabla u|^p~\!dz - (\nu_1+\nu_2)\\ &=S_{b,\gamma}
\int\limits_{{\mathcal K}}|y|^{a-p-b+\gamma}|z|^{-\gamma}~\!u^p~\!dz - (\nu_1+\nu_2)~\!,
\end{aligned}
$$
as $u$ achieves $S_{b,\gamma}$. This implies that $\nu_1=\nu_2=0$, which is equivalent to say that 
the sets $W_1, W_2$ are empty, as claimed.
\QED

\section{Proofs of the main results}
\label{S:proofs}

\subsection{Proof of Theorem \ref{T:main_H}}
\label{SS:main_H}

We start by showing that $\gamma\ge b$ is a necessary condition to have that 
$S_{b,\gamma}>0$.

Take a function $u\in {\cal 
C}^\infty_c(\R^d\setminus\Sigma_0)$ and a point 
$z_0\in\Sigma_0\cap \S^{d-1}$.  For any 
integer $h> 1$ put
$$
u_h(z)=u(z-hz_0).
$$
We estimate 
$$
S_{b,\gamma}\le \frac{\displaystyle \int\limits_{\R^d}|y|^{a}|z|^{-b}|\nabla 
u_h|^p~\!dz}
{\displaystyle \int\limits_{\R^d}|y|^{a-p-b+\gamma}|z|^{-\gamma}|u_h|^p~\!dz}
=h^{\gamma-b}~\!
\frac{\displaystyle \int\limits_{\R^d}|y|^{a}|h^{-1}z+z_0|^{-b}|\nabla 
u|^p~\!dz}
{\displaystyle 
\int\limits_{\R^d}|y|^{a-p-b+\gamma}|h^{-1}z+z_0|^{-\gamma}|u|^p~\!dz},
$$
which gives
\begin{equation}
\label{eq:first}
S_{b,\gamma} \le h^{\gamma-b}
\frac{\displaystyle \int\limits_{\R^d}|y|^{a}|\nabla u|^p~\!dz}
{\displaystyle \int\limits_{\R^d}|y|^{a-p-b+\gamma}|u|^p~\!dz}+o(h^{\gamma-b})
\end{equation}
as $h\to\infty$, because $|z_0|=1$. Thus $S_{b,\gamma}=0$ if $\gamma<b$.

To prove the converse, take $u\in {\cal C}^\infty_c(\R^d\setminus\{0\})$. We 
use 
(\ref{eq:cyl}) with $a+\gamma-b$ instead of $a$ to get 
$$
\int\limits_{\R^d}|y|^{a-p+\gamma-b}|z|^{-\gamma}|u|^p~\!dz=
\int\limits_{\R^d}|y|^{a-p+\gamma-b}\big|{|z|^{-\frac{\gamma}{p}}u}\big|^p~\!dz
\le c 
\int\limits_{\R^d}|y|^{a+\gamma-b}\big|\nabla(|z|^{-\frac{\gamma}{p}}u)\big|^p~\!dz~\!.
$$
Up to a constant, the last integral is not larger than 
$$
\int\limits_{\R^d}|y|^{a+\gamma-b}|z|^{-\gamma}\big|\nabla u\big|^p~\!dz+
\int\limits_{\R^d}|y|^{a+\gamma-b}|z|^{-\gamma-p}|u|^p~\!dz.
$$
Using (\ref{eq:CMN}) with $a+\gamma-b$ instead of $a$, we can estimate
$$
\int\limits_{\R^d}|y|^{a+\gamma-b}|z|^{-\gamma-p}|u|^p~\!dz \le c\int\limits_{\R^d}|y|^{a+\gamma-b}|z|^{-\gamma}\big|\nabla u\big|^p~\!dz,
$$
and obtain 
$$
\int\limits_{\R^d}|y|^{a-p+\gamma-b}|z|^{-\gamma}|u|^p~\!dz\le c 
\int\limits_{\R^d}|y|^{a}|z|^{-b} (|y||z|^{-1})^{\gamma-b}|\nabla u|^p~\!dz\le 
c \int\limits_{\R^d}|y|^{a}|z|^{-b} |\nabla u|^p~\!dz,
$$
because the assumption $\gamma\ge b$ trivially implies  $(|y||z|^{-1})^{\gamma-b}\le 
1$ on $\R^d\setminus\{0\}$. 
The conclusion follows, thanks to the density  Lemma \ref{L:density}.
\QED

\subsection{The infimum  $\widetilde{S_{b,\gamma}}$}
\label{S:sphere}

In this subsection we study the minimization problem (\ref{eq:min_sphere}). 

Assume $k+a>0$. Recall that $W^{1,p}(\S^{d-1};|\Pi\sigma|^ad\sigma)$ is  the 
reflexive Banach space obtained by completing ${\cal C}^\infty(\S^{d-1})$ with 
respect
to the uniformly convex norm
$$
\|\f\|_{1,p}^p=\int\limits_{\S^{d-1}}|\Pi\sigma|^{a}(|\nabla_{\!\sigma}\f|^p+|\f|^p)~\!d
\sigma.
$$

\begin{Lemma}
\label{L:main_sphere}
Assume (\ref{eq:assumptions}) and let $\widetilde{S_{b,\gamma}}$ be the infimum 
in (\ref{eq:min_sphere}).
Then the following facts hold:
\begin{itemize}
\item[$i)$] $\widetilde{S_{b,\gamma}}>0$ if and only if $\gamma\ge b$;
\item[$ii)$] If $\gamma>b$ then $\widetilde{S_{b,\gamma}}$ is achieved on 
$W^{1,p}(\S^{d-1};|\Pi\sigma|^ad\sigma)$.
\item[$iii)$] The function $\gamma\mapsto \widetilde{S_{b,\gamma}}$ is strictly 
increasing and continuous; 
\item[$iv)$] If $\gamma=b$ (which needs $k+a>p$)  then $\widetilde{S_{b,b}}\le 
\Lambda_0^p$ and the function $b\mapsto \widetilde{S_{b,b}}$ is non-increasing
and (locally) Lipschitz continuous.
\end{itemize}
\end{Lemma}

\proof
The first two statements are immediate consequences of Lemma \ref{L:compact}. 

For any fixed $\f\in W^{1,p}(\S^{d-1};\Ps^ad\sigma)\setminus\{0\}$, the functional $\widetilde{{\cal J}_{b,\gamma}}(\f)$ in (\ref{eq:J_new}) is 
evidently a strictly increasing function of $\gamma$, since $\Ps\le 1$ on $\S^{d-1}$. Thus the function 
$\gamma\mapsto \widetilde{S_{b,\gamma}}$ is non-decreasing.
For $\gamma>b$ the strict monotonicity follows as 
$\widetilde{S_{b,\gamma}}$
is achieved.

To prove the continuity of $\gamma\mapsto \widetilde{S_{b,\gamma}}$,  
fix $\gamma_0\ge b$ and $\eps>0$. Take $\f_\eps\in W^{1,p}(\S^{d-1};\Ps^ad\sigma)\setminus\{0\}$ such that $\widetilde{ {\cal J}_{b,\gamma_0}}(\f_\eps)\le \widetilde{S_{b,\gamma_0}}+\eps$.
Then  for $\gamma\searrow\gamma_0$ we have
$$
\widetilde{S_{b,\gamma_0}}\le \widetilde{S_{b,\gamma}}\le  \widetilde{ {\cal J}_{b,\gamma}}(\f_\eps)=
\widetilde{{\cal J}_{b,\gamma_0}}(\f_\eps)+o(1)\le
\widetilde{S_{b,\gamma_0}}+\eps+o(1).
$$
This proves the continuity from the right.

\medskip
Further, let $\gamma_0>b$. For any $\gamma\in(b,\gamma_0)$ let  $\f_\gamma$ be a minimizer of
$\widetilde{{\cal J}_{b,\gamma}}$. Then for $\gamma\nearrow\gamma_0$ we have
$$
\widetilde{S_{b,\gamma_0}}\ge \widetilde{S_{b,\gamma}}= \widetilde{{\cal J}_{b,\gamma}}(\f_\gamma)=
\widetilde{{\cal J}_{b,\gamma_0}}(\f_\gamma)+o(1)\ge 
\widetilde{S_{b,\gamma_0}}+o(1).
$$
This proves the continuity from the left and concludes the proof of $iii)$.

\medskip

It remains to handle the bottom case $\gamma=b$. Notice that for any given 
nontrivial $\f$, the ratio 
\begin{equation}
\label{eq:J}
\widetilde{{\cal J}_{b,b}}
(\f)=\frac{\displaystyle\int\limits_{\S^{d-1}}|\Pi\sigma|^{a}
\big[|\nabla_{\!\sigma} \f |^2+{\rmH_b^2}\f^2\big]^\frac{p}{2}~\!d\sigma}
{\displaystyle\int\limits_{\S^{d-1}}|\Pi\sigma|^{a-p}|\f |^p~\!d\sigma}
\end{equation}
decreases in $b$, because $\rmH_b=\rmH_0-\frac{b}{p}$ does. 
This 
implies that 
the function $b\mapsto  \widetilde{S_{b,b}}$ is non-increasing.

To show that  $\widetilde{S_{b,b}}\le \Lambda_0^p$ we use the trial function 
$\f_\eps(\sigma)=\Ps^{-{\Lambda_0+\eps}}$, compare with Lemma \ref{L:useful_new}. Thanks to (\ref{eq:useful})
one can easily estimate
$$
|\nabla_{\!\sigma} \f_\eps|\le 
(\Lambda_0-\eps)\Ps^{-\Lambda_0-1+\eps}+c\Ps^{-\Lambda_0+\eps}
=\Ps^{-\Lambda_0-1+\eps}\big(\Lambda_0-\eps+ c\Ps\big).
$$
Since $a-p\Lambda_0-p=-k$, we obtain that
$$
\widetilde{{\cal J}_{b,b}}(\f_\eps) \le ~
\frac{\displaystyle\int\limits_{\S^{d-1}}\Ps^{-k+p\eps} 
\big(\Lambda_0-\eps+ c\Ps\big)^p~\!d\sigma}
{\displaystyle\int\limits_{\S^{d-1}}\Ps^{-k+p\eps} ~\!d\sigma}=\Lambda_0^p+o(1)
$$
as $\eps\to 0$. Thus $\widetilde{S_{b,b}}\le \Lambda_0^p$ follows.

\medskip
Next, for fixed $\f\in {\cal C}^\infty(\S^{d-1})$ and
$\sigma\in \S^{d-1}$, the function
$$
f(\tau)= \big[|\nabla_{\!\sigma}\f|^2+\tau^2 \f^2\big]^\frac{p}{2}
$$
is non-decreasing and convex on $(0,\infty)$. Since in addition the function 
$b\mapsto \rmH_b$ decreases, 
for given  $b_2\le b_1<p\rmH_0$ we have that 
$$
f({\rmH_{b_1}})\le f({\rmH_{b_2}})\le f({\rmH_{b_1}})+ 
f'({\rmH_{b_2}})({\rmH_{b_2}}-{\rmH_{b_1}})\le 
f({\rmH_{b_1}})+\frac{p}{\rmH_{b_2}}f(\rmH_{b_2})({\rmH_{b_2}}-{\rmH_{b_1}}),
$$
because $\rmH f'({\rmH})=pf({\rmH})^{\frac{p-2}{p}}\rmH^2|\partial_rv|^2
\le pf({\rmH})$.  Since $p({\rmH_{b_2}}-{\rmH_{b_1}})=b_1-b_2$, we  infer that
$$
\widetilde{{\cal J}_{b_2,b_2}}(\f) 
\Big(1-\frac{1}{\rmH_{b_2}}~\!(b_1-
b_2)\Big)\le \widetilde{{\cal J}_{b_1,b_1}}(\f)\le \widetilde{{\cal 
J}_{b_2,b_2}}(\f).
$$
Since $\f$ was arbitrarily chosen, we get
\begin{equation}
\label{eq:continuity2}
\widetilde{S_{b_2,b_2}}-\frac{\widetilde{S_{b_2,b_2}}}{\rmH_{b_2}}~\!(b_1-
b_2)\le \widetilde{S_{b_1,b_1}}\le \widetilde{S_{b_2,b_2}}.
\end{equation}
The Lipschitz continuity of $b\mapsto \widetilde{S_{b,b}}$ is 
proved.
\QED

In the bottom case $\gamma=b$ the minimization problem (\ref{eq:min_sphere})
is non-compact, in the sense that if $\widetilde{S_{b,b}}=\Lambda_0^p$, then there 
exist 
minimizing sequences which converge
weakly to $0$. 

However, compactness can be recovered if the strict inequality 
$\widetilde{S_{b,b}}<\Lambda_0^p$ holds. This is reminiscent of
a phenomenon frequently observed  in non-compact {\it anticoercive} 
variational problems, see  the large collection of examples in the 
seminal papers \cite{Ls_comp, Ls} by P.L. Lions. 
For this reason, we state the next result in a separate lemma, the proof of which has been inspired by 
\cite{Naz, Naz2}.

\begin{Lemma} 
\label{L:widetilde_achieved}
Let $k+a>p$. If $\widetilde{S_{b,b}}<\Lambda_0^p$ then $\widetilde{S_{b,b}}$ is 
achieved.
\end{Lemma}

\proof
We first reduce the minimization problem (\ref{eq:min_sphere}) to a
1D problem. We change the coordinates $\sigma\in\S^{d-1}$ via
$$
\sigma= (t,tY,\sqrt{1-t^2}~\!X)~,\qquad t=\Ps~,\quad 
Y\in \S^{k-1}~,~~X\in\S^{d-k-1}
$$
(the variable $X$ has to be neglected if $d=k+1$).
Notice that if 
$\f$ depends only on $t$ then
\begin{equation}
\label{eq:1D}
\widetilde{{\cal J}_{b,b}}(\f)=\widehat{ {\cal J}_{b,b}}(\f):=\frac{\displaystyle\int\limits_0^1 
{t}^{k+a-1}(1-{t}^2)^{\frac{d-k-2}2}\big[(1-{t}^2)|\f'|^2+{\rmH_b^2}\f^2\big]^\frac{p}{2}~\!d{t}}
{\displaystyle\int\limits_0^1 {t}^{k+a-1-p}(1-{t}^2)^{\frac{d-k-2}2}|\f|^p~\!d{t}}.
\end{equation}
Denote by
$W^{1,p}_{\rm 1D}(\S^{d-1};\Ps^ad\sigma)$ the space of functions in 
$W^{1,p}(\S^{d-1};\Ps^ad\sigma)$ depending only on the $t$ variable.
We claim that in fact
\begin{equation}
\label{eq:hat}
\widetilde{{S}_{b,b}}=\widehat{{S}_{b,b}}:=\inf_{\f\in W_{\rm 1D}^{1,p}(\S^{d-1};\Ps^ad\sigma)} \widehat{ {\cal J}_{b,b}}(\f).
\end{equation}
Indeed, for arbitrary  $\f\in  W^{1,p}(\S^{d-1};\Ps^ad\sigma)$ and for 
fixed 
$(Y, X)\in \S^{k-1}\times\S^{d-k-1}$ we have that 
$$
\int\limits_0^1 {t}^{k+a-1}(1-{t}^2)^{\frac{d-k-2}2}\big[(1-{t}^2)|\partial_t\f|^2+{\rmH_b^2}\f^2\big]^\frac{p}{2}~\!d{t}
\ge \widehat{S_{b,b}} \int\limits_0^1 
{t}^{k+a-1-p}(1-{t}^2)^{\frac{d-k-2}2}|\f|^p~\!d{t}.
$$
Since  $|\nabla_{\!\sigma}\f|\ge |\partial_t\f|$, integration in the remaining 
variables $(Y,X)$ gives
$$
\int\limits_{\S^{d-1}}|\Pi\sigma|^{a}\big[|\nabla_{\!\sigma} \f |^2+{\rmH_b^2}\f^2\big]^\frac{p}{2}~\!d\sigma \ge \widehat{S_{b,b}} 
\int\limits_{\S^{d-1}}|\Pi\sigma|^{a-p}|\f |^p~\!d\sigma.
$$
Thus $\widetilde{S_{b,b}}\ge \widehat{S_{b,b}}$. The 
opposite inequality is trivial, and the claim follows.

\medskip

The remaining part of the proof is based on the idea of \cite{Ls_comp, Ls}, see 
also \cite[Lemma 3.1]{Naz}. We take a minimizing sequence 
$\f_j\in W_{\rm 1D}^{1,p}(\S^{d-1};\Ps^ad\sigma)$. Moreover, we can suppose that the following 
properties hold:

\begin{description}
\item[$(a)$] $\f_j\ge0$; 
 \item[$(b)$] 
$\displaystyle{\int\limits_0^1 
{t}^{k+a-1-p}(1-{t}^2)^{\frac{d-k-2}2}|\f_j|^p~\!d{t}} 
=1$,\\ so that
$
\displaystyle{\int\limits_0^1 
{t}^{k+a-1}(1-{t}^2)^{\frac{d-k-2}2}\big((1-{t}^2)|\f_j'|^2+{\rmH_b^2}\f_j^2\big
)^\frac{p}{2}~\!d{t} }=\widetilde{S_{b,b}}+o(1);
$

\item[$(c)$] $\f_j \to \bm\f$  weakly in $W_{\rm 
1D}^{1,p}(\S^{d-1};|\Pi\sigma|^ad\sigma)$ and in 
$L^p(\S^{d-1};|\Pi\sigma|^{a-p}d\sigma)$;

\item[$(d)$] 
The sequences
\begin{equation}
\label{eq:measures}
{t}^{k+a-1-p}(1-{t}^2)^{\frac{d-k-2}2}|\f_j|^p~,\qquad  
{t}^{k+a-1}(1-{t}^2)^{\frac{d-k-2}2}\big((1-{t}^2)|\f_j'|^2+{\rmH_b^2}\f_j^2\big
)^\frac{p}{2}
\end{equation}
converge weakly in the space of measures on $[0,1]$.
\end{description}

\noindent
Since $\f_j \to {\bm\f}$ strongly in 
$L^p_{\mathrm{loc}}(\S^{d-1}\setminus\Sigma_0;\Ps^{a-p})$, we have 
\begin{equation}
\label{eq:Lions}
{t}^{k+a-1-p}(1-{t}^2)^{\frac{d-k-2}2}|\f_j|^p~\!dt~\rightharpoondown ~
{t}^{k+a-1-p}(1-{t}^2)^{\frac{d-k-2}2}|{\bm\f}|^p~\!dt+A{\delta_0}(t)
\end{equation}
for some $A\ge 0$, where ${\delta_0}$ is the Dirac delta.
Notice that $(b)$ implies
\begin{equation}
\label{eq:denom}
\int\limits_0^1{t}^{k+a-1-p}(1-{t}^2)^{\frac{d-k-2}2}|{\bm\f}|^p~\!d{t}+A=1.
\end{equation}

Next, let $\mathfrak{m}$ be the limit measure relative to the second sequence 
in 
(\ref{eq:measures}).
For any fixed function
$\phi\in {\cal C}^\infty_c(\S^{d-1}\setminus\Sigma_0)$ depending only on $t$ we 
obtain
$\phi\nabla_{\!\sigma} \f_j \to \phi\nabla_{\!\sigma} {\bm\f}$ weakly in 
$L^p(\S^{d-1};|\Pi\sigma|^ad\sigma)$, and thus
$$
\aligned
\int\limits_0^1|\phi|^pd\mathfrak{m}= &\,
\lim_{j\to\infty}\int\limits_0^1 {t}^{k+a-1}(1-{t}^2)^{\frac{d-k-2}2}\big[(1-{t}^2)|\f'_j|^2+{\rmH_b^2}\f_j^2\big]^\frac{p}{2}|\phi|^p~\!d{t}
\\
\ge&\, \int\limits_0^1 {t}^{k+a-1}(1-{t}^2)^{\frac{d-k-2}2}\big[(1-{t}^2)|{\bm\f}'|^2+{\rmH_b^2}{\bm\f}^2\big]^\frac{p}{2}|\phi|^p~\!
d{t}
\endaligned
$$
by the convexity of the functional involved. That is 
\begin{equation}
\label{eq:mu_1} 
\mathfrak{m}\ge{t}^{k+a-1}(1-{t}^2)^{\frac{d-k-2}2}\big[(1-{t}^2)|{\bm\f}'|^2+{\rmH_b^2}{\bm\f}^2\big]^\frac{p}{2}~\!. 
\end{equation}

Now let $\psi=\psi(t)$ be a smooth cutoff function such that
$$
\psi\equiv 1\ \ \mbox{in} \ \ [0,\eps]; \quad \psi\equiv 0\ \ \mbox{in} \ \ 
[2\eps,1]; \quad |\psi'|\le c\eps^{-1}
$$
and recall the elementary inequality
$$
\aligned
\big[(1-{t}^2){\f'_j}{}^2+{\rmH_b^2}\f_j^2\big]^\frac{p}{2}~\!|\psi|^p\ge 
&\,(1-{t}^2)^\frac { p } { 2 } |
\f'_j\psi|^p\\
\ge &\, (1-{t}^2)^\frac{p}{2}\Big[|(\f_j \psi)'|^p
-c(p)\left(|\f_j\psi'|^p+|\f'_j\psi|^{p-1}~\!|\f_j\psi'|\right)\Big].
\endaligned
$$
We infer that
\begin{equation}
\label{eq:leading}
\aligned
\int\limits_0^1|\psi|^pd\mathfrak{m}= &\,
\lim_{j\to\infty}\int\limits_0^1 {t}^{k+a-1}(1-{t}^2)^{\frac{d-k-2}2}\big[(1-{t}^2)|\f'_j|^2+{\rmH_b^2}\f_j^2\big]^\frac{p}{2}|\psi|^p~\!d{t}
\\
\ge&\, (1-(2\eps)^2)^{\frac{d-k-2+p}2}\lim_{j\to\infty}\int\limits_0^1 {t}^{k+a-1}|(\f_j\psi)'|^p
~\!d{t}
-c(p)\,\lim_{j\to\infty}\left(I_{1j}+I_{2j}\right),
\endaligned
\end{equation}
where 
$$
I_{1j}=\int\limits_0^1 {t}^{k+a-1}|\f_j\psi'|^p~\!d{t}~; \qquad
I_{2j}= \int\limits_0^1 {t}^{k+a-1}|\f'_j\psi|^{p-1}~\!|\f_j\psi'|~\!d{t}.
$$
We estimate the main term in (\ref{eq:leading}) via the 1D-Hardy inequality, which gives
$$
\int\limits_0^1 {t}^{k+a-1}|(\f_j\psi)'|^p
~\!d{t} \ge \Lambda_0^p\int\limits_0^1 {t}^{k+a-1-p}|\f_j\psi|^p~\!d{t}.
$$
Therefore, from (\ref{eq:Lions}) we infer that
$$
\lim_{j\to\infty}\int\limits_0^1 {t}^{k+a-1}|(\f_j\psi)'|^p
~\!d{t}\ge \Lambda_0^p A.
$$
Next, we use H\"older inequality to estimate
$$
I_{2j}\le \Big(\int\limits_0^1 {t}^{k+a-1}|\f'_j\psi|^p~\!d{t}\Big)^{\frac {p-1}p}\cdot I_{1j}^{\frac 1p}\le cI_{1j}^{\frac 1p},
$$
where $c$ depends only on the 
sequence 
$\f_j$.

It remains to estimate $I_{1j}$. Since $\f_j \to {\bm\f}$ strongly in 
$L^p(\S^{d-1};\Ps^a)$ by Lemma \ref{L:compact}, 
we have that 
$$
\lim_{j\to\infty}I_{1j}=\int\limits_0^1 
{t}^{k+a-1}|{\bm\f}\psi'|^p
~\!d{t}\le c\int\limits_\eps^{2\eps} 
{t}^{k+a-1}|\eps^{-1}{\bm\f}|^p
~\!d{t}\le c\int\limits_\eps^{2\eps} 
{t}^{k+a-1-p}|{\bm\f}|^p
~\!d{t}.
$$
The latter quantity vanishes  as $\eps\to0$. Therefore, 
$\lim\limits_{j\to\infty} I_{1j}=\lim\limits_{j\to\infty} I_{2j}=o(1)$ as $\eps\to 0$.

Summing up, we get
$$
\int\limits_0^1|\psi|^pd\mathfrak{m}\ge 
(1-(2\eps)^2)^{\frac{d-k-2+p}2}\Lambda_0^p 
A-o(1),
$$
as $\eps\to 0$, 
that gives $\mathfrak{m}\ge \Lambda_0^p A{\delta_0}$. Comparing this inequality 
with (\ref{eq:mu_1}) we obtain
$$ 
\mathfrak{m}\ge{t}^{k+a-1}(1-{t}^2)^{\frac{d-k-2}2}\big[(1-{t}^2)|{\bm\f}'|^2+{\rmH_b^2}{\bm\f}^2\big]^\frac{p}{2}+\Lambda_0^p 
A{\delta_0}. 
$$
Since $\f_j$ is minimizing, we have
$$
\aligned
\widetilde{S_{b,b}}=\int\limits_0^1d\mathfrak{m} &\, \ge 
\int\limits_0^1 {t}^{k+a-1}(1-{t}^2)^{\frac{d-k-2}2}\big[(1-{t}^2)|{\bm\f}'|^2+{\rmH_b^2}{\bm\f}^2\big]^\frac{p}{2}~\!d{t}+
\Lambda_0^p A\\
&\, \ge 
\widetilde{S_{b,b}}\int\limits_0^1{t}^{k+a-1-p}(1-{t}^2)^{\frac{d-k-2}2}|{\bm\f}|^p~\!d{ t}+\Lambda_0^pA=
\widetilde{S_{b,b}}(1-A)+\Lambda_0^pA,
\endaligned
$$
see (\ref{eq:denom}). 
Since $\widetilde{S_{b,b}}<\Lambda_0^p$, we 
infer 
that $A=0$ and that $\bm\f$ is a minimizer for the functional 
$\widehat{{\cal J}_{b,b}}$, which ends the proof.
\QED

\subsection{Proof of Theorem \ref{T:Main2}}
\label{SS:Main2}

It suffices to prove that $S_{b,\gamma}=\widetilde{S_{b,\gamma}}$, as the continuity and the monotonicity of $\gamma\mapsto S_{b,\gamma}$
then follows by $iii)$ in Lemma \ref{L:main_sphere}.

We start by using $ii)$ in Lemma \ref{L:main_sphere} to find a minimizer $\Phi\in W^{1,p}(\S^{d-1};\Ps^ad\sigma)$ for $\widetilde{S_{b,\gamma}}$.
Then we test $S_{b,\gamma}$ with the function defined in spherical coordinates by
$$
u_\delta(r\sigma)=\begin{cases}
r^{-\rmH_b+ \delta}\Phi(\sigma)&
\text{if $r\le 1$,}\\
r^{-\rmH_b- \delta}\Phi(\sigma)&
\text{if $r> 1$,}
\end{cases}
$$
for any $\delta>0$. Easily, $u_\delta\in \mathcal D^{1,p}(\R^d;|y|^a|z|^{-b} dz)$. By adapting the direct computations in
the proof of \cite[Theorem 3]{CMN}, one infers that
$$
\begin{aligned}
\int\limits_{\R^d}|y|^{a-p-b+\gamma}|z|^{-\gamma}| u_\delta|^p~\!dz&=
\frac{2}{p\delta}
~\!\int\limits_{\mathbb S^{d-1}}|\Pi\sigma|^{a-p- b +\gamma}|\Phi|^p~\!d\sigma~\!,\\
\int\limits_{\R^d}|y|^a|z|^{-b}|\nabla u_\delta|^p~\!dz&=\frac{2}{p\delta}
\Big(~\!\int\limits_{\mathbb S^{d-1}}|\Pi\sigma|^{a}\big[|\nabla_{\!\sigma}\Phi|^2+\rmH_b^2\Phi^2\big]^\frac{p}{2}~\!d\sigma+O(\delta^2)\Big)
\end{aligned}
$$
as $\delta\to 0$.
Therefore
$$
S_{b,\gamma}\le \frac{\displaystyle\int\limits_{\R^d}|y|^a|z|^{-b}|\nabla u_\delta|^p~\!dz}
{\displaystyle\int\limits_{\R^d}|y|^{a-p-b+\gamma}|z|^{-\gamma}| u_\delta|^p~\!dz}=   \widetilde{S_{b,\gamma}}+O(\delta^2).
$$
This proves that $S_{b,\gamma}\le \widetilde{S_{b,\gamma}}$.

To prove the opposite inequality use a standard argument to show   that $\Phi$ can not change sign in ${{\S^{d-1}}}$. Thus, we can assume 
that 
$\Phi$ is nonnegative, and solves the 
corresponding Euler--Lagrange equation in a weak sense.

Consider the function defined in spherical coordinates by
$$
U(r\sigma)=r^{-\rmH_b}\Phi(\sigma).
$$
Evidently $U\in W^{1,p}_{\rm loc}(\mathcal K;|y|^a|z|^{-b}dz)$, where $\mathcal K$ is the cone in (\ref{eq:cone}). 
Arguing as in \cite[Theorem 3]{CMN}, one can prove that $U$
is a local weak solution of
$$
-\div(|y|^a|z|^{-b}|\nabla U|^{p-2}\nabla U)= 
\widetilde{S_{b,\gamma}}|y|^{a-p-b+\gamma}|z|^{-\gamma}|U|^{p-2}U\qquad \text{in 
$\R^d\setminus\{0\}$.}
$$
By Lemma \ref{L:lsc}, $U$ is lower semicontinuous and positive on $\mathcal K$.
Thus Theorem \ref{T:AP} applies  and gives $S_{b,\gamma}\ge \widetilde{S_{b,\gamma}}$, which ends the proof.
\QED

\subsection{Proof of Theorem \ref{T:Main3}}
\label{SS:Main3}

We use the linear transform $u\mapsto \widetilde{u}$ defined in spherical 
coordinates by
\begin{equation}
\label{eq:ttau}
u(r\sigma)=\widetilde{u}(r^{\tau_b}\sigma)~,\quad \text{where 
$\tau_b:=\frac{\rmH_b}{\rmH_0}$,}
\end{equation}
see for instance \cite{Ho}.
Let $u\in {\cal C}^\infty_c(\R^d\setminus\Sigma_0)$. We compute 
\begin{gather}
\label{eq:Ttau_a}
\int\limits_{\R^d}|y|^{a}|z|^{-b}|\nabla u|^p~\!dz= \tfrac{\rmH_0}{\rmH_b}
\int\limits_{0}^\infty r^{d+a-1}dr \int\limits_{\S^{d-1}} 
\Ps^{a}\Big[\tfrac{\rmH_b^2}{\rmH_0^2}|\partial_r \widetilde{u}|^2+ 
\!r^{-2}|\nabla_{\!\sigma} \widetilde{u}|^2\Big]^\frac{p}{2} d\sigma~\!,
\\
\nonumber
\int\limits_{\R^d}|y|^{a-p}|z|^{-b}|u|^p~\!dz=\tfrac{\rmH_0}{\rmH_b}\int\limits_
{0}^\infty r^{d+a-p-1}dr \int\limits_{\S^{d-1}}  \Ps^{a-p}|\widetilde{u}|^p ~\!d\sigma~\!.
\end{gather}

Notice that for $b=0$ the integrals in (\ref{eq:Ttau_a}) provide the $p$-th 
power of the norm in $\mathcal D^{1,p}(\R^d;|y|^adz)$
(purely cylindrical weights). If $b\neq 0$ then the right-hand side of 
(\ref{eq:Ttau_a}) gives the $p$-th power of an equivalent norm
in  $\mathcal D^{1,p}(\R^d;|y|^adz)$. Therefore, the transform $u\mapsto 
\widetilde{u}$ can be uniquely extended to an invertible isomorphism
$$
\mathcal D^{1,p}(\R^d;|y|^a|z|^{-b}dz)\longrightarrow \mathcal 
D^{1,p}(\R^d;|y|^a dz).
$$
In addition, we infer that
\begin{equation}
\label{eq:equivalent}
S_{b,b}=\inf_{v\in D^{1,p}(\R^d;|y|^a dz) } \mathcal R_{\tau_b}(v)
\end{equation}
where, for $\tau>0$, we put
$$
\mathcal R_{\tau}(v)= \frac{\displaystyle\int\limits_{0}^\infty 
r^{d+a-1}dr\int\limits_{\S^{d-1}} \Ps^{a}
\Big[\tau^2 |\partial_r v|^2+ \!r^{-2}|\nabla_{\!\sigma} v|^2\Big]^\frac{p}{2} 
d\sigma}
{\displaystyle\int\limits_{0}^\infty r^{d+a-p-1}dr\int\limits_{\S^{d-1}} 
\Ps^{a-p}|v|^p d\sigma}~\!.
$$
Next, for fixed $v\in {\cal C}^\infty_c(\R^d)$, $r>0$ and 
$\sigma\in \S^{d-1}$ the function
$f(\tau)= \big[\tau^2 |\partial_r v|^2+ \!r^{-2}|\nabla_{\!\sigma} v|^2\Big]^\frac{p}{2}
$
is non-decreasing and convex on $(0,\infty)$. By repeating the argument we 
used 
to prove (\ref{eq:continuity2}) one can show that for $b_2\le b_1<p\rmH_0$ it holds that
$$
S_{b_2,b_2}-\frac{S_{b_2,b_2}}{\rmH_{b_2}}~\!(b_1-b_2)\le S_{b_1,b_1}\le 
S_{b_2,b_2}~\!,
$$
which ends the proof of $i)$ in Theorem \ref{T:Main3}.

\medskip

Before proving $ii)$ we recall that the space ${\cal C}^\infty_c(\R^d\setminus\Sigma_0)$ is dense
in $\mathcal D^{1,p}(\R^d;|y|^a|z|^{-b}dz)$ and in $\mathcal D^{1,p}(\R^d;|y|^a dz)$, use
Lemma \ref{L:density}.
Therefore,  formula (\ref{eq:first}) with $\gamma=b$ and (\ref{eq:cyl})
imply
$$
S_{b,b}=S_{b,b}(\R^d\setminus\Sigma_0)\le  \inf_{u\in {\cal C}^\infty_c(\R^d\setminus\Sigma_0)}
\frac{\displaystyle~ \int\limits_{\R^d}|y|^{a}|\nabla u|^p~\!dz}
{\displaystyle \int\limits_{\R^d}|y|^{a-p}|u|^p~\!dz}= S_{0,0}=\Lambda_0^p~\!.
$$
Since $S_{b,b}$ is non-increasing with $b$, we infer that 
$S_{b,b}=\Lambda_0^p$ for any $b\le 0$.
On the other hand, 
$
\lim\limits_{b\to p\rmH_0} S_{b,b}=0
$
because for $b=p\rmH_0$ the weight in the numerator of the ratio in 
(\ref{eq:Sbb}) is still locally integrable, while the weight in the 
denominator is not. Therefore,
the conclusion of the proof of $ii)$
follows by the continuity and the monotonicity of the function $b\mapsto S_{b,b}$.

Now we prove the equality $S_{b,b} =\widetilde{S_{b,b}}$. 
For any $\delta>0$ we have that 
$$
S_{b,b}\le S_{b,b+\delta} =\widetilde{S_{b,b+\delta}},
$$
use Theorem \ref{T:Main2}.
Therefore, for any given   $\f\in W^{1,p}(\S^{d-1};\Ps^a 
d\sigma)\setminus\{0\}$ 
we have 
$$
S_{b,b}
\le \frac{~\displaystyle\int\limits_{\S^{d-1}}|\Pi\sigma|^{a}\big[|\nabla_{\!\sigma} \f |^2+{\rmH_b^2}\f^2\big]^\frac{p}{2}~\!d\sigma}
{\displaystyle\int\limits_{\S^{d-1}}|\Pi\sigma|^{a-p+\delta}|\f |^p~\!d\sigma}=\frac{~\displaystyle\int\limits_{\S^{d-1}}|\Pi\sigma|^{a}
\big[|\nabla_{\!\sigma} \f |^2+{\rmH_b^2}\f^2\big]^\frac{p}{2}~\!d\sigma}
{\displaystyle\int\limits_{\S^{d-1}}|\Pi\sigma|^{a-p}|\f |^p~\!d\sigma}+o(1)
$$
as $\delta\to 0$.
This proves that $S_{b,b}\le \widetilde{S_{b,b}}$ for any $b<p\rmH_0$.

\medskip

Next, we just proved that $S_{b,b}=\Lambda_0^p$ for $b\le b_*$. By $iv)$ in Lemma 
\ref{L:main_sphere} we have that 
$\widetilde{S_{b,b}}\le \Lambda_0^p$ for any $b$. Hence $S_{b,b}=\widetilde{S_{b,b}}=\Lambda_0^p$ for $b\le 
b_*$.

\medskip

By testing $\widetilde{S_{b,b}}$ with the constant function $\f\equiv 1$ one 
infers that $\lim\limits_{b\nearrow 
p\rmH_0}\widetilde{S_{b,b}}=\lim\limits_{\rmH_b\searrow 
0}\widetilde{S_{b,b}}=0$.
Thus  $\widetilde{S_{b,b}}<\Lambda_0^p$ for $b$ close enough to $p\rmH_0$. Put
$$
\beta_*=\inf\{b<p\rmH_0~|~ \widetilde{S_{b,b}}<\Lambda_0^p\}\ge b_*.
$$
In fact, we have that $\beta_*=b_*$. Indeed, for any $b\in (\beta_*,p\rmH_0)$ 
it 
turns out that $\widetilde{S_{b,b}}$ is achieved by Lemma 
\ref{L:widetilde_achieved}.
By repeating literally the first part of the proof of Theorem \ref{T:Main2},
one gets that $S_{b,b}=\widetilde{S_{b,b}}<\Lambda_0^p$. Since $b>\beta_*$ was 
arbitrarily chosen, 
by continuity we have 
$S_{\beta_*,\beta_*}=\widetilde{S_{\beta_*,\beta_*}}=\Lambda_0^p$. Thus 
$\beta_*= b_*$.
This ends the proof of the equality $S_{b,b} =\widetilde{S_{b,b}}$. 

Finally, recall that
for any fixed $\f\in W^{1,p}(\S^{d-1};\Ps^ad\sigma)\setminus\{0\}$, the function $b\mapsto \widetilde{{\cal J}_{b,b}}(\f)$ in (\ref{eq:J}) is 
strictly decreasing. Since in addition $\widetilde{S_{b,b}}$ is achieved for $b\in(b_*,p\rmH_0)$,
we infer that the map $b\mapsto \widetilde{S_{b,b}}={S_{b,b}}$ is strictly decreasing on $[b_*,p\rmH_0)$,
 which concludes the proof.
\QED

\subsection{Proof of Theorem \ref{T:non-achiev}}
\label{SS:na}

We start with the case $\gamma=b<b_*$.  
Recall that the transform (\ref{eq:ttau}) 
gives ${\cal J}_{b,b}(u)=\mathcal R_{\tau_b}(\widetilde u)$ and that the minimization problems in (\ref{eq:Sbb}) and in 
(\ref{eq:equivalent}) are equivalent. 

Since $\tau_b=\rmH_b/\rmH_0$ strictly decreases as $b$ increases, we infer that
the function 
$b\mapsto \mathcal R_{\tau_b}(v)$ is strictly decreasing for any fixed $v\in \mathcal D^{1,p}(\R^d;|y|^a dz)\setminus\{0\}$. 
If $S_{b_0,b_0}$ were achieved for some $b_0<b_*$, then $\mathcal R_{\tau_{b_0}}$ would also attain its minimum, and the 
function $b\mapsto S_{b,b}$ would be strictly increasing in a right  neighborhood of 
$b_0$. But $S_{b,b}=\Lambda_0^p$ for 
any  $b<b_*$ by $ii)$ in Theorem \ref{T:Main3}, a contradiction.

Next, let $\gamma>b$ or $\gamma=b>b_*$. Then the infimum 
$\widetilde{S_{b,\gamma}}$ is achieved, see $ii)$ in Lemma 
\ref{L:main_sphere} and Lemma \ref{L:widetilde_achieved}, respectively. Without 
loss of generality, the minimizer 
$\Phi\in W^{1,p}(\S^{d-1};\Ps^ad\sigma)$ 
is nonnegative. Since 
$\widetilde{S_{b,\gamma}}=S_{b,\gamma}$, as in the proof 
of Theorem \ref{T:Main2}  one can show that the function 
$U(r\sigma)=r^{-\rmH_b}\Phi(\sigma)$ belongs to $W^{1,p}_{\rm loc}(\mathcal K;|y|^a|z|^{-b}dz)$ (the cone $\mathcal K$ is introduced in (\ref{eq:cone}))
and 
solves the equation (\ref{eq:eq_mS}) in $\R^d\setminus\{0\}$. Since $U\notin 
 D^{1,p}(\R^d;|y|^a|z|^{-b}dz)$, 
thanks to Lemma \ref{L:AP2} we infer  
that $S_{b,\gamma}$ is not achieved in this case.

It remains to consider the case $\gamma=b=b_*$. 

Similarly to the proof of $ii)~\Rightarrow~i)$ in Theorem \ref{T:AP} (case $k+a\ge p$),
one can find a lower 
semicontinuous
and positive weak solution $\Phi\in W^{1,p}_{\rm loc}({\S^{d-1}\setminus\Sigma_0};\Ps^ad\sigma)$ of the 
Euler--Lagrange equation corresponding to the functional $\widetilde{{\cal J}_{b_*,b_*}}$. As before, this implies that the function 
$U(r\sigma)=r^{-\rmH_{b_*}}\Phi(\sigma)$
solves the equation (\ref{eq:eq_mS}) in $\R^d\setminus\Sigma_0$. 
Since $U$ satisfies the assumptions in Lemma \ref{L:AP2} we infer  
that $S_{b_*,b_*}$ is not achieved. 
\QED

\begin{Remark}
\label{R:final_tilde}
 As for the attainability of $\widetilde{S_{b,\gamma}}$, we know 
the following facts:
\begin{itemize}
\item[$i)$]  If $\gamma>b$ then $\widetilde{S_{b,\gamma}}$ is achieved by $ii)$ 
in Lemma \ref{L:main_sphere};
\item[$ii)$] If $b\in(b_*,p\rmH_0)$ then 
$\widetilde{S_{b,b}}$ is achieved by Lemma \ref{L:widetilde_achieved};
\item[$iii)$] If $b<b_*$ then $\widetilde{S_{b,b}}$ and is 
not achieved (argue as for $S_{b,b}$).
\end{itemize}

We do not know whether $\widetilde{S_{b_*,b_*}}$ is achieved, except for the 
case $p=2$. In the latter case the function $\Phi_0(\sigma)=\Ps^{-\Lambda_0}$
solves the corresponding Euler--Lagrange equation in $\S^{d-1}\setminus\Sigma_0$, cf. Theorem \ref{T:sharp} below.
Since $\Phi_0\notin W^{1,2}(\S^{d-1}; \Ps^ad\sigma)$, by arguing as in 
Lemma \ref{L:AP2} one can infer that $\widetilde{S_{b_*,b_*}}$ is not achieved.

By arguing as in Theorem \ref{T:AP}, for general   $p>1$ one can prove the existence of a nonnegative solution 
$\Phi\in W^{1,p}_{\rm loc}({\S^{d-1}\setminus\Sigma_0})\setminus\{0\}$
of the corresponding
Euler--Lagrange equation. We conjecture that 
$\Phi\notin W^{1,p}({\S^{d-1}\setminus\Sigma_0},\Ps^ad\sigma)$,
which would imply that $\widetilde{S_{b_*,b_*}}$ is not 
achieved. 
\end{Remark}

\section{Some estimates on $b_*$}
\label{S:b*}

In this section we assume that $k+a>p$, so that we can deal with the {\it bottom case} $\gamma=b$. The parameter
$b_*\in[0,p\rmH_0)$ is introduced in Theorem \ref{T:Main3}.

\begin{Theorem}
\label{T:b*positive}
Assume $d>k+1$ or $p>\frac32$. Then $b_*>0$.
\end{Theorem}

\proof
For $\alpha\ge 0$ to be determined later, we define
$$
U(z)=|y|^{-{\Lambda_0}}|z|^{-\alpha{\Lambda_0}}.
$$
Our goal is to show that if $b>0$ is small enough one can choose 
$\alpha$ in such a 
way that 
\begin{equation}
\label{eq:inequality}
\mathcal LU:=-\div(|y|^a|z|^{-b}|\nabla U|^{p-2}\nabla U)- {\Lambda_0^p} 
|y|^{a-p}|z|^{-b}U^{p-1}\ge 0\quad \text{on $\R^d\setminus\Sigma_0$},
\end{equation}
(and in fact equality holds if 
$b=\alpha=0$). Then Theorem \ref{T:AP} 
 implies $S_{b,b}\ge\Lambda_0^p$ for $b>0$ small, 
 and part $ii)$ in Theorem \ref{T:Main3}
gives the conclusion of the proof.

\medskip

It is convenient to put 
$$
t=t(z)=|y||z|^{-1}\in(0,1]~, 
\qquad F_\alpha=F_\alpha(t)= (1+(\alpha^2+2\alpha)t^2)^\frac{p-2}{2}
$$
where, with some abuse of notation, we identify $y$ with $(y,0)=\Pi 
z$. We easily have 
$y\cdot z=|y|^2=|z|^2 t^2$ and
\begin{equation}
\label{eq:yz}
\nabla(|y|^A|z|^B)\cdot y= |y|^A|z|^B\big(A+B 
t^2)~,\qquad
\nabla(|y|^A|z|^B)\cdot z= |y|^A|z|^B\big(A+B)
\end{equation}
for any $A, B\in\R$. Also we keep in mind that 
\begin{equation}
\label{eq:U^p-1}
|y|^{\Lambda_0-k}|z|^{-b-\alpha{\Lambda_0}(p-1)}=|y|^{a-p}|z|^{-b}U^{p-1}.
\end{equation}

We calculate
$$
\nabla U=-\Lambda_0|y|^{-1}U\big(|y|^{-1}y+\alpha t|z|^{-1}z\big)\ \ 
\Longrightarrow \ \
-|\nabla U|^{p-2}\nabla 
U={\Lambda_0^{p-1}}|y|^{1-p}U^{p-1}F_\alpha\big(|y|^{-1}y+\alpha 
t|z|^{-1}z\big),
$$
and thus
$$
-\div(|y|^a|z|^{-b}|\nabla U|^{p-2}\nabla 
U)={\Lambda_0^{p-1}}(\mathcal A_1+\alpha~\!\mathcal A_2),
$$
where
$$
\mathcal A_1=\div(|y|^{a-p}|z|^{-b}U^{p-1}F_\alpha y)~,\quad 
\mathcal A_2=\div(|y|^{a-p+2}|z|^{-b-2}U^{p-1}F_\alpha z).
$$
To compute $\mathcal A_1$ we first use the first equality in (\ref{eq:yz}) with 
$A=a-p-{\Lambda_0}(p-1)={\Lambda_0}-k$ and $B=-b-{\Lambda_0}(p-1)\alpha$. 
Taking into account (\ref{eq:U^p-1}) we get
$$
\nabla\big(|y|^{a-p}|z|^{-b}U^{p-1}\big)\cdot y= 
|y|^{a-p}|z|^{-b}U^{p-1}\big({\Lambda_0}-k-(b+{\Lambda_0}(p-1)\alpha)t^2\big).
$$
Also, (\ref{eq:yz}) gives
$$
\nabla F_\alpha\cdot 
y=\frac{p-2}{2}F_\alpha^{\frac{p-4}{p-2}}(\alpha^2+2\alpha)\nabla\big(|y|^2|z|^{
-2} \big)\cdot y
=(p-2) F_\alpha^{\frac{p-4}{p-2}}(\alpha^2+2\alpha)(1-t^2)t^2,
$$
and we infer that
$$
\begin{aligned}
\mathcal A_1&=F_\alpha \nabla\big(|y|^{a-p}|z|^{-b}U^{p-1}\big)\cdot y+ 
|y|^{a-p}|z|^{-b}U^{p-1}\big(\nabla F_\alpha\cdot y +k F_\alpha\big)\\
&= |y|^{a-p}|z|^{-b}U^{p-1}\Big({\Lambda_0} F_\alpha-(b+{\Lambda_0} (p-1)\alpha)F_\alpha t^2+
(p-2)F_\alpha^{\frac{p-4}{p-2}}(\alpha^2+2\alpha)(1-t^2)t^2\Big).
\end{aligned}
$$

To compute $\mathcal A_2$,
we use the second identity in (\ref{eq:yz}) and (\ref{eq:U^p-1}) to get
$$
\nabla\big(|y|^{a-p+2}|z|^{-b-2}U^{p-1}\big)\cdot z= 
|y|^{a-p}|z|^{-b}U^{p-1}\big({\Lambda_0}-k-(b+{\Lambda_0}(p-1)\alpha)\big)t^2.
$$
Since (\ref{eq:yz}) for $A=-B=2$ gives also $\nabla F_\alpha\cdot z=0$, we
finally we arrive at
$$
\begin{aligned}
\mathcal A_2&=F_\alpha \nabla\big(|y|^{a-p+2}|z|^{-b-2}U^{p-1}\big)\cdot z+ 
d|y|^{a-p+2}|z|^{-b-2}U^{p-1} F_\alpha\\
&= |y|^{a-p}|z|^{-b}U^{p-1}
({\Lambda_0}+d-k-b-{\Lambda_0}(p-1)\alpha)F_\alpha t^2.
\end{aligned}
$$
Therefore,
\begin{equation}
\aligned
\label{eq:end0}
\mathcal LU=&\,{\Lambda_0^{p-1}}\big(\mathcal A_1+\alpha\mathcal 
A_2-\Lambda_0|y|^{a-p}|z|^{-b}U^{p-1}\big)={\Lambda_0^{p-1}}
|y|^{a-p}|z|^{-b}U^{p-1}\times\\
\times&\,\Big[{\Lambda_0}\big(F_\alpha-1-(p-2)\alpha 
F_\alpha t^2)
+\big(C_\alpha F_\alpha +(p-2)(\alpha^2+2\alpha)F_\alpha^{\frac{p-4}{p-2}}~(1-t^2)\big)t^2 \Big],
\endaligned
\end{equation}
where
$C_\alpha=(d-k-b)\alpha-b-{\Lambda_0}(p-1)\alpha^2$.

Now we put $b=\eps\alpha$, where $\eps$ is a constant to be chosen 
later. This gives $C_\alpha=(d-k-\eps)\alpha+O(\alpha^2)$.
Since for $t\in(0,1]$ and small $\alpha$ we trivially 
have $F_\alpha=1+(p-2)\alpha t^2+O(\alpha^2t^2)$, the expression is square 
brackets in (\ref{eq:end0}) can be rewritten as follows:
\begin{equation}
\label{eq:end2}
\alpha t^2\Big[(d-k-\eps)+2(p-2)~(1-t^2)+O(\alpha)\Big].
\end{equation}
It remains to take a small $\eps\in(0,1)$ such that
$d-k-\eps+2(p-2)>0$ (recall that we are assuming $p>3/2$ if $d-k=1$). Then for 
any $\alpha>0$ small enough (hence, for any $b>0$ small enough) the expression 
in (\ref{eq:end2}) is nonnegative for all $t\in(0,1]$. This gives 
(\ref{eq:inequality}) and completes the proof.
\QED

\begin{Remark}
We conjecture that $b_*$ is always strictly positive, i.e., that the restriction $p>\frac32$ in the case $d=k+1$ purely technical.
\end{Remark}

We can also provide some (not sharp) estimates from above on $b_*$.

\begin{Proposition} It holds that 
$b^*\le p(\rmH_0-\Lambda_0)=d-k$.
\end{Proposition}

\proof
We test $\widetilde{S_{b,b}}$ with the function $\f(\sigma)=\Ps^{-\Lambda_0+\delta}$, which is indeed a function of $t=\Ps\in [0,1]$. 
Using formulae (\ref{eq:1D}) and (\ref{eq:hat})
we have that
{\begin{equation}
\label{**}
\widetilde{S_{b,b}}\le 
\frac{~\displaystyle\int\limits_0^1 
{{t}}^{-1+p\delta}(1-{{t}}^2)^{\frac{d-k-2}2}\big[
(\Lambda_0-\delta)^2(1-{{t}}^2)+{\rmH_b^2} {{t}}^2\big] ^\frac{p}{2}~\!d{{t}}}
{\displaystyle\int\limits_0^1 
{{t}}^{-1+p\delta}(1-{{t}}^2)^{\frac{d-k-2}2}~\!d{{t}}}=:s_\delta(b).
\end{equation}
Since $p\rmH_b=p\rmH_0-b$, it is evident from (\ref{**}) that the  the function 
$s_\delta(\cdot)$ is decreasing 
with respect to $b$. For $b=(\rmH_0-(\Lambda_0-\delta))p$
we have $\rmH_b=\Lambda_0-\delta$, hence (\ref{**}) gives  $s_\delta(b)=(\Lambda_0-\delta)^p$. This means that for 
$b>(\rmH_0-\Lambda_0)p$ and small $\delta$ we have that $s_\delta(b)<\Lambda_0^p$, which concludes the proof.
}
\QED

In the Hilbertian case $p=2$ we have a complete picture.

\begin{Theorem}
\label{T:sharp}
Let $p=2$. Then $b_*$ is the smallest root of the equation
$b^2-4\rmH_0 b+4(\rmH_0-\Lambda_0)^2=0$, that is,
\begin{equation}
\label{eq:b_star}
b_*=2\rmH_0-2\sqrt{\rmH_0^2-(\rmH_0-\Lambda_0)^2}
\in (0,2\rmH_0).
\end{equation}
Moreover, for $b\in (b_*,2\rmH_0)$ it holds that 
\begin{equation}
\label{eq:sharp}
S_{b,b}=
\Lambda_0^{2}-\big(\sqrt{\rmH_0^2-\rmH_b^2}-(\rmH_0-\Lambda_0)\big)^2.
\end{equation}
\end{Theorem}

\proof 
Recall that in the Hilbertian case, we have
$2\rmH_b= (d+a-2)-b=2\rmH_0-b$, $2\Lambda_0={k+a-2}$.

For any $b<2\rmH_0$, $\lambda\in\R$ the function 
$$
U_{b,\lambda}(z)=|y|^{-\lambda}|z|^{\lambda-\rmH_b}
$$
is smooth on  $\R^d\setminus\Sigma_0$. Computations similar to those in the proof of Theorem \ref{T:b*positive} show that 
$U_{b,\lambda}$ solves
\begin{equation}
\label{eq:p=2}
-\div\big(|y|^a|z|^{-b}\nabla 
U\big)=\lambda(2\Lambda_0-\lambda)|y|^{a-2}|z|^{-b}U
+\big[(\rmH_b-\lambda)^2-b\lambda\big]|y|^{a}|z|^{-b-2}U.
\end{equation}

Let $\mathfrak b$ be the constant in the right hand side of (\ref{eq:b_star}). We choose $\lambda=\Lambda_0$. 
Since $(\rmH_{\mathfrak b}-\Lambda_0)^2=\mathfrak b \Lambda_0$,  from (\ref{eq:p=2}) we see that the function
$U_{\mathfrak b,\Lambda_0}$ solves
$$
-\div\big(|y|^a|z|^{-\mathfrak b}\nabla 
U\big)=\Lambda_0^{2}|y|^{a-2}|z|^{-\mathfrak b}U\qquad \text{in 
$\R^d\setminus\Sigma_0$}.
$$
Then 
Theorem \ref{T:AP} and Corollary \ref{C:equal} give $S_{\mathfrak b,\mathfrak b} 
\ge \Lambda_0^{2}$, hence 
$S_{\mathfrak b,\mathfrak b}= \Lambda_0^{2}$ and $\mathfrak b\le b_*$ by Theorem 
\ref{T:Main3}.

\medskip

If $b\in(\mathfrak b,2\rmH_0)$ we denote by $\mathfrak S_{b}
$ the constant in the right 
hand 
side of (\ref{eq:sharp}) and put 
\begin{equation}
\label{eq:lambda_b}
\lambda_b:=\rmH_0-\sqrt{\rmH_0^2-\rmH_b^2}\in(0,\Lambda_0).
\end{equation}
Notice that $\lambda_b^2-2\rmH_0\lambda_b+\rmH_b^2=0$ and
$$
(\rmH_b-\lambda_b)^2=b\lambda_b~,\qquad 
\mathfrak S_{b}=\Lambda_0^{2}-(\Lambda_0-\lambda_b)^2=\lambda_b(2\Lambda_0-\lambda_b).
$$
Therefore, from (\ref{eq:p=2}) we see that  $U_{b,\lambda_b}$ solves
$$
-\div\big(|y|^a|z|^{- b}\nabla U\big)=\mathfrak S_{b}|y|^{a-2}|z|^{-b}U
\quad \text{in $\R^d\setminus\Sigma_0$,}
$$
which implies $S_{b,b} \ge \mathfrak S_{b}$, thanks to Theorem \ref{T:AP} and 
Corollary \ref{C:equal}.

\medskip

Recall that $S_{b,b}=\widetilde{S_{b,b}}$ by Theorem \ref{T:Main3}. Therefore, to
conclude the proof it suffices to show that $\widetilde{S_{b,b}}\le \mathfrak S_{b}$. We 
test $\widetilde{S_{b,b}}$ with the function 
$$
\f_b(\sigma)=\Ps^{-\lambda_b},
$$
where $\lambda_b<\Lambda_0$ is defined in (\ref{eq:lambda_b}) (notice that $\f_b\in W^{1,2}(\S^{d-1};\Ps^ad\sigma)$ by Lemma \ref{L:useful_new}). We use (\ref{eq:useful}) to compute
$$
\begin{aligned}
\int\limits_{\S^{d-1}}|\Pi\sigma|^{a}\big[|\nabla_{\!\sigma} \f_b |^2+\rmH_b^2|\f_b|^2\big]~\!d\sigma&=
\int\limits_{\S^{d-1}}|\Pi\sigma|^{a-2\lambda_b-2}\big[\lambda_b^2|\nabla_{\!\sigma} \Ps|^2+\rmH_b^2\Ps^{2}\big]~\!d\sigma
\\
&=\lambda_b^2\int\limits_{\S^{d-1}}|\Pi\sigma|^{a-2\lambda_b-2}~\!d\sigma
+(\rmH_b^2-\lambda_b^2)\int\limits_{\S^{d-1}}|\Pi\sigma|^{a-2\lambda_b}~\!d\sigma~\!.
\end{aligned}
$$
Then, Lemma \ref{L:useful_iii} gives
$$
\int\limits_{\S^{d-1}}|\Pi\sigma|^{a}\big[|\nabla_{\!\sigma} \f_b |^2+\rmH_b^2|\f_b|^2\big]~\!d\sigma=\Big(\lambda_b^2 + (\rmH_b^2-\lambda_b^2)\frac{\Lambda_0-\lambda_b}{\rmH_0-\lambda_b}\Big)
\int\limits_{\S^{d-1}}|\Pi\sigma|^{a-2-2\lambda_b}~\!d\sigma~\!. 
$$
We now use again the identity  $\lambda_b^2-2\rmH_0\lambda_b+\rmH_b^2=0$ to write 
$\rmH_b^2-\lambda_b^2=2\lambda_b(\rmH_0-\lambda_b)$. We infer that the constant in front of the last integral  is simply 
$\mathfrak S_b$. In conclusion, we have that 
$$
\int\limits_{\S^{d-1}}|\Pi\sigma|^{a}\big[|\nabla_{\!\sigma} \f_b |^2+\rmH_b^2|\f_b|^2\big]~\!d\sigma=
\mathfrak S_b\int\limits_{\S^{d-1}}|\Pi\sigma|^{a-2}|\f_b|^2~\!d\sigma~\!.
$$
We infer that $\mathfrak S_b\ge \widetilde{S_{b,b}}$, hence equality holds. 
\QED

\begin{Remark}
%\label{R:p=2}
We indeed proved that for $b\in (b_*,2\rmH_0)$ the function $\f_b$ is a minimizer for 
$\widetilde{S_{b,b}}$.

Notice that for  $p=2$ we know the explicit value of $S_{b,\gamma}$ if $\gamma=b+p$ and if $\gamma=b$, see \cite{CMN} and  Theorem \ref{T:sharp}.
The remaining cases are open.

\end{Remark}

\end{document}